 \documentstyle [12pt,twoside, epsf]{article}
 \newcommand{\ZZ}{\mbox{$Z\!\!\! Z\!$}} 

 \let\\=\cr
 
 \def\Chi{\hbox{\raise0.5ex\hbox{$\chi$}}}

 \newtheorem{th}{Theorem}
 \newtheorem{lem}{Lemma}
 \newtheorem{cor}{Corollary}

 \newtheorem{defn}{Definition}

 \newtheorem{rem}{Remark}

 \def\picill#1by#2(#3){\epsffile{#3}}

\begin{document}
 \pagestyle{myheadings}

  \markboth{{\sc Kauffman \& Lambropoulou}}{{\sc On the classification of
 rational knots}}

 \title{On the Classification of Rational Knots}

 \author{Louis H. Kauffman and Sofia Lambropoulou}

\date{}

 \maketitle

\begin{abstract} 
\noindent In this paper we give combinatorial proofs of the classification of  unoriented
and oriented rational knots  based on the  now known classification of alternating knots 
and the calculus of continued fractions. We also characterize the class of
strongly invertible rational links.  Rational links are of fundamental importance in
the study of DNA recombination. 
\end{abstract}

\noindent {\it AMS Subject Classification: \ 57M27}

\noindent {\it Keywords:} \ rational knots and links, rational tangles, 
continued fractions, tangle fraction, flypes,  alternating knots and links, chirality,
invertibility.

 \section{Introduction}

Rational knots and links comprise the simplest class of links.  The first
twenty five knots, except for $8_5$, are rational.  Furthermore  all knots and links up to ten
crossings are either rational or are obtained by inserting rational tangles into a small
number of planar graphs, see \cite{C1}. Rational links are alternating
with one or two unknotted components,  and they are also known in the literature as
Viergeflechte, four-plats  or $2$-bridge knots depending on their geometric
representation. More precisely, rational knots can be represented as: 
\smallbreak

-  plat closures of four-strand braids (Viergeflechte \cite{BS}, four-plats). These are knot
diagrams with two local maxima and two local minima. 
\smallbreak

-  $2$-bridge knots. A $2$-bridge knot is a knot that has
a diagram in which there are two distinct arcs, each overpassing a
consecutive sequence of crossings, and every crossing in the diagram is in
one of these sequences. The two arcs are called the bridges of the
diagram (compare with \cite{BZ}, p. 23).
\smallbreak

-  numerator or denominator closures of rational tangles (see Figures 1, 5).  A rational tangle is
the result of consecutive twists on  neighbouring endpoints of two trivial arcs.  For examples see
Figure 1 and Figure 3.
\smallbreak

All three representations are equivalent. The equivalence between the first and the
third is easy to see by planar isotopies. For the equivalence between the first and the second
representation see for example \cite{BZ}, pp. 23, 24.  In this paper we consider
rational  knots as obtained by taking numerator or denominator closures of rational tangles (see
Figure 5).

\bigbreak

$$ \picill4.3inby2.2in(R1)  $$

\begin{center}
{Figure 1 - A rational tangle and a rational knot } 
\end{center}
\vspace{3mm}


The notion of a tangle was introduced in 1967 by Conway  \cite{C1} in his work on enumerating
and classifying knots and links, and he defined the rational knots as numerator or denominator
closures  of the rational tangles. (It is worth noting here that 
 Figure 2 in \cite{BS} illustrates a rational tangle, but no
special importance is given to this object. It is obtained from a four-strand braid by
plat-closing only the top four ends.)  Conway \cite{C1}  also defined  {\it the fraction}
of a rational tangle to be a rational number or $\infty.$ He observed that this
 number for a rational tangle equals a continued fraction expression with all
numerators equal to one and all denominators of the same sign, that can be read from a
tangle diagram in alternating standard form. Rational tangles are classified by their
fractions by means of the following theorem. 

\begin{th}[Conway, 1975]\label{Conway}{ \ Two rational tangles are isotopic if and
only if they have the same fraction.
  } \end{th}

Proofs of Theorem 1 are given in \cite{Mo}, \cite{BZ} p.196, \cite{GK2} and \cite{KL1}. 
The first two proofs  invoked the classification of rational knots and the theory of
branched covering spaces. The $2$-fold branched
covering spaces of $S^3$ along the rational links give rise to the lens spaces $L(p,q)$,
see \cite{Sei}. The proof in \cite{GK2} is the first combinatorial  proof  of
this theorem. The proofs in \cite{Mo}, \cite{BZ} and \cite{GK2} use definitions different
from the above for the fraction of a rational tangle. In \cite{KL1} a new combinatorial
proof of Theorem~1  is given using the solution of the  Tait
Conjecture for alternating knots  \cite{Ta}, \cite{MT} adapted for tangles. 
A second combinatorial proof  is
given in \cite{KL1} using coloring for defining the tangle fraction. 
\smallbreak

 Throughout the paper by the term `knots' we will refer to both
knots and links, and whenever we really mean `knot' we shall emphasize it. 
More than one rational tangle can yield the same or isotopic  
rational knots and the equivalence relation between the rational tangles is mapped  
into an arithmetic equivalence of their corresponding fractions. Indeed we have the
following.

\begin{th}[Schubert, 1956]\label{Schubert1}{ \ Suppose that rational tangles with
fractions $\frac{p}{q}$ and $\frac{p'}{q'}$ are given ($p$ and $q$ are relatively prime.
Similarly for $p'$ and $q'$.) If  $K(\frac{p}{q})$ and $K(\frac{p'}{q'})$
denote the corresponding rational knots obtained by taking numerator closures of
these tangles, then $K(\frac{p}{q})$ and $K(\frac{p'}{q'})$ are isotopic 
if and only if

\begin{enumerate}
\item $p=p'$ and
\item either $q\equiv q'\, mod \, p$  \ or \ $qq'\equiv 1\, mod \, p.$
\end{enumerate}
} \end{th}

Schubert \cite{Sch2} originally stated the  classification of rational knots  and links 
by representing them  as $2$-bridge links. Theorem \ref{Schubert1} has hitherto been
proved by taking  the $2$-fold branched covering spaces of $S^3$ along $2$-bridge links,
showing that these correspond bijectively to oriented diffeomorphism classes of lens
spaces, and invoking the classification of lens spaces \cite{Rd3}. 
Another proof  using covering spaces has been given by  Burde in \cite{Bu}. See also the
excellent notes  on the subject by Siebenmann \cite{Sie}.
 The above statement of Schubert's theorem is a formulation of the Theorem in the
language of Conway's tangles.

\smallbreak
Using his methods for the unoriented case, Schubert also extended the classification of
rational knots  and links to the case of oriented rational knots and links described as
$2$-bridge links.  Here is our formulation of the Oriented Schubert Theorem written in
the language of Conway's tangles. 

\begin{th}[Schubert, 1956]\label{Schubert2}{ \ Suppose that orientation-compatible
rational tangles with fractions $\frac{p}{q}$  and $\frac{p'}{q'}$ are given with $q$ and
$q'$ odd. ($p$ and $q$ are relatively prime. Similarly for $p'$ and $q'$.) If
$K(\frac{p}{q})$ and $K(\frac{p'}{q'})$ denote the corresponding rational knots obtained
by taking numerator closures of these tangles, then $K(\frac{p}{q})$ and
$K(\frac{p'}{q'})$ are isotopic if and only if
\begin{enumerate}
\item $p=p'$ and
\item either $q\equiv q'\, mod \, 2p$  \ or \ $qq'\equiv 1\, mod \, 2p.$
\end{enumerate}
} \end{th}

Theorems \ref{Schubert1} and \ref{Schubert2} could have been stated equivalently using the
denominator closures of rational tangles. Then the arithmetic
equivalences of the tangle fractions related to isotopic knots would be the same as in
Theorems
\ref{Schubert1} and \ref{Schubert2}, but with the roles of numerators and denominators
exchanged. 

\smallbreak

This paper  gives the first combinatorial  proofs of  Theorems
\ref{Schubert1} and \ref{Schubert2} using tangle theory. Our proof of Theorem
\ref{Schubert1} uses the results and the techniques developed in \cite{KL1}, while the
proof of  Theorem
\ref{Schubert2} is based on that of Theorem \ref{Schubert1}.  We have
located the essential points in the proof of the classification of rational knots in the 
 question: {\it Which rational tangles will close
 to form a specific knot or link diagram?} By looking at the Theorems in
 this way, we obtain a path to the results that can be understood without
 extensive background in three-dimensional topology. In the course of these proofs we see
connections between the elementary number theory of fractions and continued fractions,
and the topology of knots and links. In order to compose these proofs  we  use the fact
that rational knots are alternating (which follows from the fact that rational tangles
are alternating, and for which we believe we found the simplest possible proof, see
\cite{KL1}, Proposition 2). We then rely on the {\it Tait Conjecture} \cite{Ta}  concerning the
classification of alternating knots, which  states the following:  
\smallbreak
\noindent {\it Two alternating knots 
 are isotopic if and only if any two corresponding reduced diagrams on $S^2$ are  related by a
finite sequence of flypes (see Figure 6).}
\smallbreak

 \noindent A diagram is said to be {\it reduced} if at every crossing the
four local regions indicated at the crossing are actually parts of four distinct global  regions in
the diagram (See \cite{Li} p. 42.). It is not hard to see that any knot or link has reduced
diagrams that represent its isotopy class. The conjecture was posed by P.G. Tait, \cite{Ta} in
1877 and was proved   by W.~Menasco and M.~Thistlethwaite, \cite{MT} in 1993. Tait did not
actually  phrase this statement as a conjecture. It was a working hypothesis for his efforts in
classifying knots.
\smallbreak

Our proof of the Schubert Theorem is elementary upon assuming
the Tait Conjecture, but this is easily stated and understood.  This paper
will be  of interest to mathematicians and biologists.

\bigbreak
The paper is organized as follows. In Section 2 we give the general set up for rational
tangles, their isotopies and operations, as well as their association to a
continued fraction isotopy invariant. In this section we also recall the
basic theory  and a canonical form of continued fractions. 
In Section 3 we prove Theorem
\ref{Schubert1} about  the classification of unoriented rational knots by means of a
direct  combinatorial and arithmetical analysis  of rational knot diagrams, using the
classification of rational tangles and the Tait Conjecture.  
In Section 4  we discuss
chirality of knots and give a classification of  the achiral rational knots and links as
numerator closures of even palindromic rational tangles  in continued fraction form
(Theorem \ref{mirror}).
In Section 5  we discuss the connectivity patterns of the four end arcs of rational
tangles and we relate connectivity to the parity of the fraction of a rational tangle
(Theorem \ref{connectivity}).
 In Section 6  we give our interpretation of the statement  of
Theorem \ref{Schubert2} and  we prove  the classification of oriented rational knots,
using the  methods we developed in the unoriented case and examining the connectivity
patterns of oriented rational  knots. In Section 6  it is
pointed out that all  oriented rational knots and  links are invertible (reverse the
orientation of both components).  
In  Section 7 
we give a classification of the strongly invertible rational links (reverse the
orientation of one component) as closures of odd  palindromic oriented rational tangles
in continued fraction form (Theorem \ref{strongly}).

\bigbreak
Here is a short history of the theory of rational knots. As explained in  
\cite{Gr}, rational knots and links were first considered by O. Simony  in 1882,
\cite{Si1, Si2, Si3, Si4}, taking twistings and knottings of a band. Simony \cite{Si2} was  
the first one to relate knots to continued fractions. After about sixty years Tietze
wrote  a series of papers \cite{Ti1, Ti2, Ti3, Ti4}  with reference
to Simony's work. Reidemeister \cite{Rd2} in 1929 calculated 
the knot group of a special class of four-plats (Viergeflechte), but four-plats were
really studied  by Goeritz \cite{Goe} and by Bankwitz and Schumann \cite{BS} in 1934. In
\cite{Goe} and \cite{BS}  proofs are given independently and with different techniques
that  rational knots have $3$-strand-braid representations, in the sense that the first
strand of the four-strand braids can be free of crossings, and that they are
alternating. (See Figure 20 for an example and Figure 26 for an abstract $3$-strand-braid
representation.)  The proof of the latter in
\cite{BS} can be easily  applied on the corresponding rational tangles in
standard form. (See Figure 1 for an example and Figure 8 for abstract representations.) 
\smallbreak

 In  1954 Schubert \cite{Sch1}  introduced the bridge representation of knots. He then
showed that the four-plats are exactly the knots that can be represented by diagrams
with  two bridges and consequently he classified rational knots by finding canonical forms
via representing them as $2$-bridge knots, see\cite{Sch2}. His proof was based on
Seifert's observation that the $2$-fold branched coverings of $2$-bridge knots \cite{Sei}
give rise  to lens spaces and on the classification of lens spaces by Reidemeister
\cite{Rd3} using Reidemeister torsion and following the lead of \cite{ST} 
(and later by Brody \cite{Br} using the knot theory of the lens space).  See also \cite{PY}.)  
Rational knots and rational tangles figure prominently in the applications of knot theory to the
topology of DNA, see \cite{Su}. Treatments of various  aspects of rational knots and rational
tangles can be found  in many places in the literature, see for example \cite{C1}, \cite{Sie},
\cite{R},  \cite{BZ}, \cite{BM},  \cite{M},  \cite{Kaw}, \cite{Li}.

\section{Rational Tangles and their Invariant Fractions}

In this section we recall from \cite{KL1} the facts that we need about rational tangles,
continued fractions and the classification of rational tangles. We intend the paper to be 
as  self-contained as possible. 
\smallbreak
 A {\it $2$-tangle} is a proper embedding  of
two unoriented arcs  and a finite number of circles in a $3$-ball
$B^3,$ so that the four endpoints lie in the boundary of $B^3.$ 
A   {\it rational tangle}  is a proper embedding  of two unoriented arcs  $\alpha_1,
\alpha_2$  in a $3$-ball $B^3,$ so that the four endpoints lie in the boundary of $B^3,$ 
and such that there exists a homeomorphism of pairs: 

$$\overline{h}: (B^3, \alpha_1, \alpha_2)  \longrightarrow (D^2\times I, \, \{x,y\}\times I) {\mbox
\ \ (a \ trivial \ tangle).}$$
  
\noindent This is equivalent to saying that rational  tangles have specific
representatives obtained by applying a finite number of consecutive twists of
neighbouring endpoints  starting from two  unknotted and unlinked arcs. Such a pair of
arcs comprise the $[0]$ or $[\infty]$ tangles, depending on their position in the plane,
see illustrations in Figure 2. We shall use this characterizing  property of a rational
tangle as our definition, and we shall then say  that the rational tangle is in {\it
twist form}. See  Figure 3 for an example.

\bigbreak

To see the equivalence of the above definitions, let $S^2$
denote the two-dimensional sphere, which is the boundary of the $3$-ball $B^3$  and let $p$ denote 
four specified points  in
$S^2.$ Let further \ $h: (S^2,p) \longrightarrow (S^2,p)$ \ be a self-homeomorphism of
$S^2$  with the four points. This extends to a self-homeomorphism
$\overline{h}$ of the $3$-ball $B^3$ (see \cite{R}, page 10).  Further, let $a$ denote the
two straight arcs $\{x,y\}\times I$  joining pairs of the four points in the boundary of
$B^3.$ Consider now $\overline{h}(a).$ We call this the tangle induced by $h.$ We note that
up  to isotopy (see definition below) $h$ is a composition of braidings of pairs of
points in $S^2$ (see \cite{PS}, pages 61 to 65).  Each such braiding induces a twist in
the corresponding tangle. So, if $h$ is a composition of braidings of pairs of points,
then the extension $\overline{h}$ is a composition of twists of neighbouring end arcs. 
Thus $\overline{h}(a)$ is a rational tangle and every rational tangle can be obtained this
way.  
\bigbreak

$$ \picill5.3inby1.5in(R2) $$

\begin{center}
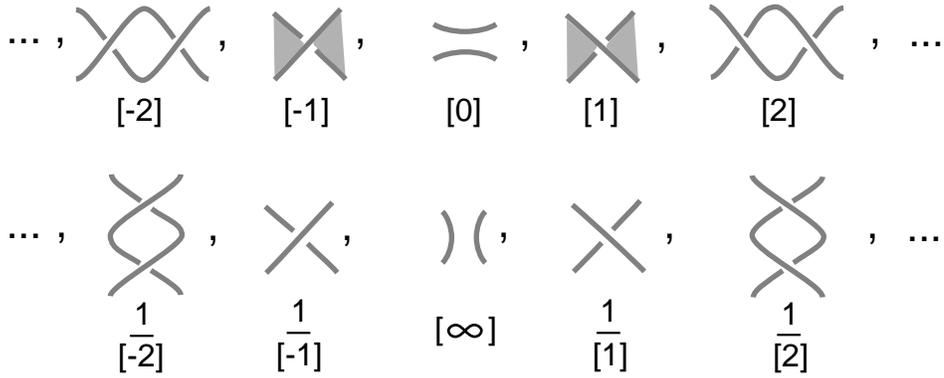
{Figure 2 - The elementary rational tangles and the types of crossings } 
\end{center}
\vspace{3mm} 


\bigbreak

$$ \picill3.55inby3.6in(R3) $$

\begin{center}
{Figure 3 - A rational tangle in twist form } 
\end{center}
\vspace{3mm} 


A {\it tangle diagram} is a regular projection of the tangle on a meridinal disc. 
 Throughout the paper by `tangle'  we will mean `regular tangle diagram'. 
 The type of crossings of knots and $2$-tangles follow
 the checkerboard  rule: shade the regions of the tangle (knot) in two
 colors, starting from the left (outside) to the right (inside) with grey, and so that
adjacent regions have different colors. Crossings in the tangle
 are said to be of {\it positive type} if they are arranged with respect to the shading 
as exemplified in Figure 2 by the tangle $[+1],$ i.e. they have the  region on the right
shaded as one walks towards the crossing along the over-arc. Crossings of the reverse
type are said to be of {\it negative type} and they are exemplified in Figure 2 by the
tangle $[-1].$ The reader should note that our crossing type  conventions are
the opposite of those of Conway in \cite{C1} and of those of Kawauchi in \cite{Kaw}.  Our
conventions agree with those of  Ernst and Sumners \cite{ES3}, \cite{Su} which in turn follow the
standard conventions of  biologists.

\smallbreak We  are
interested in tangles up to isotopy.  Two rational tangles, $T, S$, in $B^3$ are {\it
isotopic}, denoted by $T \sim S$,  if and only if any two  diagrams of
them have identical configurations of their four endpoints on the boundary of the
projection disc, and they differ by a finite sequence of the well-known Reidemeister moves
\cite{Rd2}, which take place in the interior of the  disc.  Of course, each  twisting
operation used in the definition of a rational tangle changes the isotopy class of  the tangle to which
it is applied.

\bigbreak
\noindent {\bf $2$-Tangle operations.} The symmetry of the four endpoints of $2$-tangles
allows for the following well-defined (up to isotopy) operations in the class of
$2$-tangles, as described in Figure 4. We have {\it the sum} of two $2$-tangles,
denoted by `$+$' and {\it the product} of two
$2$-tangles, denoted  by `$*$'. This product `$*$' is not to be confused with Conway's 
product `$\cdot$' in \cite{C1}. 

In view of these operations we can say that a rational
tangle is created inductively by consecutive additions of the tangles
$[\pm 1]$ on the right or on the left and multiplications  by the  tangles $[\pm 1]$ at
the bottom or at the top, starting from the tangles $[0]$ or $[\infty].$ And since, when
we start creating a rational tangle, the very first  crossing  can be
 equally seen as a horizontal or as a vertical one, we may always assume 
that we start twisting from the tangle $[0].$  
 Addition and multiplication of tangles are
not commutative. Also, they do not preserve  the class of rational tangles. The sum
(product) of two rational tangles is rational if and only if one of the two consists in a
number of horizontal (vertical) twists.   
\smallbreak
 The {\it mirror image} of a tangle $T,$  denoted $-T,$ is obtained from
$T$ by switching  all the crossings. So we have $-[n] = [-n]$ and 
$ -\frac{1}{[n]} =\frac{1}{[-n]}.$  Finally, {\it the rotation} of $T$, denoted 
$T^r$, is obtained by rotating $T$ on its plane counterclockwise by $90^0,$ whilst   {\it the
inverse} of $T$, denoted $T^i$,   is defined to be  ${-T}^r$.  Thus inversion is
 accomplished by rotation and mirror image.  For example, ${[n]}^i = \frac{1}{[n]} $ and
${\frac{1}{[n]}}^i = {[n]}.$  Note  that $T^r$ and $T^i$
are in general not isotopic to $T$. 

\bigbreak

 $$\vbox{\picill4.95inby1.8in(R4)  }$$

\begin{center}
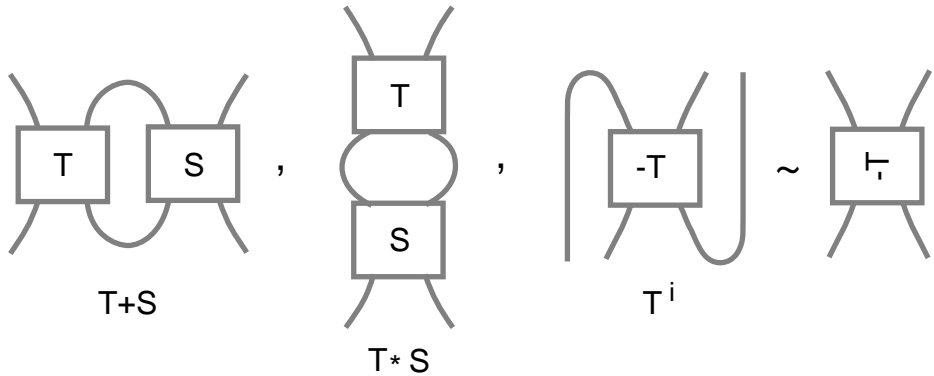
{Figure 4 - Addition, multiplication and inversion of $2$-tangles } 
\end{center}
\vspace{3mm} 


  Moreover, by joining with simple arcs
the two upper and the two lower endpoints of a $2$-tangle $T,$ we obtain a 
knot called the {\it Numerator} of $T$, denoted by $N(T).$ Joining with simple arcs  each
pair of the corresponding top and bottom endpoints of $T$ we obtain the {\it Denominator}
of $T,$  denoted by $D(T)$. We have
$N(T) = D(T^r)$ and $D(T) = N(T^r).$  We point out that
the numerator closure of the sum of two rational tangles  is
still a rational knot or link. But the denominator closure of the sum of
two rational tangles  is not necessarily a rational knot or link, think for example of the sum 
$\frac{1}{[3]} + \frac{1}{[3]}.$

\bigbreak

 $$\vbox{\picill5.45inby1.4in(R5)  }$$

\begin{center}
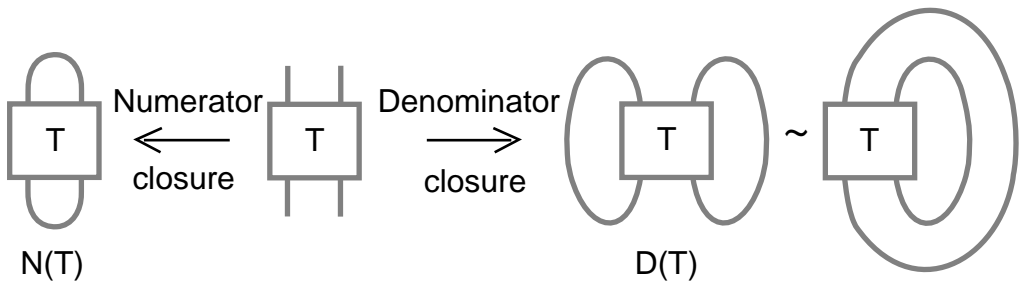
{Figure 5 - The numerator and denominator of a $2$-tangle  } 
\end{center}
\vspace{3mm}


\noindent {\bf Rational tangle isotopies.} We define now two isotopy moves for rational
tangles that play a crucial role in the theory of rational knots and rational tangles. 

\begin{defn}{\rm \  A {\it flype} is an isotopy of a  $2$-tangle $T$ (or a knot or link) applied on a
$2$-subtangle of the form $[\pm 1]+t$ or $[\pm 1]*t$ as illustrated in Figure 6. A flype
fixes the endpoints of the subtangle on which it is applied.  A flype shall be called
{\it rational} if the $2$-subtangle on which it applies is rational. }
\end{defn}

\bigbreak

$$ \picill3.65inby2.1in(R6) $$

\begin{center}
{Figure 6 - The flype moves  } 
\end{center}
\vspace{3mm} 


 We define the {\it truncation} of a rational tangle to be the result of
partially untwisting the tangle. For rational tangles,
flypes are of very specific types. Indeed, we have:
\smallbreak

\noindent {\it Let  $T$ be a  rational tangle in twist form. Then} 
{\it \begin{itemize}
\item[(i)] \ $T$ does not contain any non-rational $2$-subtangles.  
\vspace{-.1in}
\item[(ii)] \ Every $2$-subtangle of $T$ is a truncation of $T$.
\end{itemize} }

\noindent For a proof of these statements we
refer the reader to our paper \cite{KL1}. As a corollary we have that all flypes of a
rational tangle $T$ are rational. 

\begin{defn}{\rm \  A {\it flip} is a rotation in space  of a $2$-tangle by $180^0$. We
say that $T^{hflip}$ is the {\it horizontal flip} of the $2$-tangle $T$ if $T^{hflip}$ is
obtained from $T$ by a $180^0$ rotation around a horizontal axis on the plane of $T$, and
$T^{vflip}$ is the {\it vertical flip} of the tangle $T$ if $T^{vflip}$ is obtained from
$T$ by a $180^0$ rotation around a vertical axis on the plane of $T.$ See Figure 7 for
illustrations.   } \end{defn}  

\bigbreak

$$ \picill3.67inby2.4in(R7) $$

\begin{center}
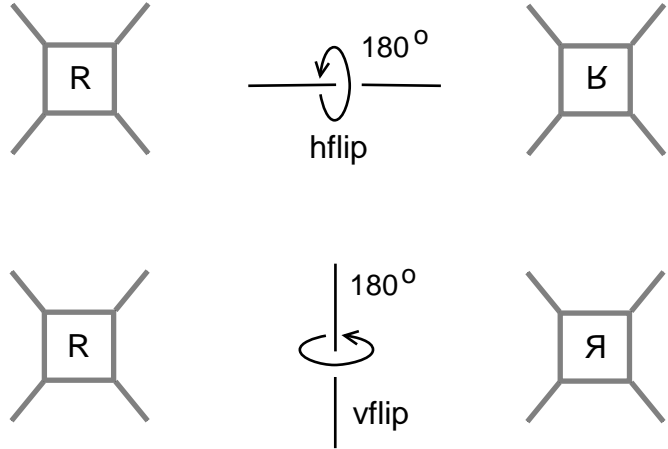
{Figure 7 - The horizontal and the vertical flip   } 
\end{center}
\vspace{3mm}


 \noindent  Note that a flip switches the endpoints of the tangle and, in general, a
flipped tangle is not isotopic to the original one; the following is a remarkable
property of rational tangles:  

\bigbreak
\noindent {\bf The Flipping Lemma} {\it  \  If \, $T$ is rational, then: 
\[ (i) \ T \sim T^{hflip}, \ \ \ \ (ii) \ T \sim T^{vflip} \ \ \ \ \mbox{and} \ \ \ \ (iii) \ T
\sim (T^{i})^{i} = (T^{r})^{r}.\]
 }

\noindent To see (i) and (ii) we apply induction and a sequence of flypes, see \cite{KL1}
for details.  $(T^{i})^{i} = (T^{r})^{r}$ is the tangle obtained from $T$ by rotating
it on its plane by $180^0,$ so statement (iii) follows by applying a vertical flip after a
horizontal flip. Note that the above statements are obvious for the tangles $[0], [\infty], \,
[n]$ and $\frac{1}{[n]}.$ 
 Statement (iii) says that {\it for rational  tangles the inversion is an operation of
order 2.}   For this reason we shall denote the inverse of a rational tangle $T$ by $1/T,$
and hence the rotation of the tangle $T$ will be denoted by $-1/T.$ This explains the
notation for the tangles $\frac{1}{[n]}.$  For arbitrary $2$-tangles the inversion is an
order 4 operation. Another consequence of the above property is that addition and
multiplication by  $[\pm 1]$ are commutative.
\bigbreak

\noindent {\bf Standard form, continued fraction form and canonical form for rational
tangles.}   Recall that the twists generating the rational tangles could take place
between the right, left, top or bottom endpoints of a previously created rational tangle.
Using obvious flypes on appropriate subtangles one can always bring the twists all to the
right (or all to the left) and to the bottom (or to the top) of the tangle. We shall then say that
the  rational tangle is in {\it standard form}.  For example Figure 1 illustrates the tangle 
$(([3]*\frac{1}{[-2]})+[2])$ in standard form. In order to read out the standard form of a
rational tangle in twist form we transcribe it as an algebraic sum using 
horizontal and vertical twists.  For example, Figure 3 illustrates the tangle $ ( (
([3] * \frac{1}{[3]}) + [-1]) * \frac{1}{[-4]}) + [2]$ in non-standard form.  
\bigbreak 

Figure 8 illustrates two equivalent (by the Flipping Lemma) ways of  representing an
abstract rational tangle in standard form: the  {\it standard representation} of a
rational tangle. In either illustration the rational tangle begins  to twist from the
tangle $[a_n]$ ($[a_5]$ in Figure 8), and it untwists from the tangle $[a_1].$ 
Note that the tangle in Figure 8 has an odd number of sets of twists ($n=5$) and this causes 
$[a_1]$ to be horizontal. If $n$ is
even and $[a_n]$ is horizontal then $[a_1]$ has to be vertical. 
\smallbreak

Another way of representing an abstract rational tangle in standard form is illustrated
in Figure 9. This is the {\it $3$-strand-braid representation}.  For an example
see Figure 10. As Figure 9 shows, the $3$-strand-braid representation is actually a
compressed version of the standard representation, so the two representations are
equivalent by a planar rotation. The upper row of crossings of the $3$-strand-braid  representation
corresponds to the horizontal crossings of the standard representation  and the  lower row to the
vertical ones. Note that, even though
the type of crossings does not change by this planar rotation, we need to draw the mirror
images of the even terms, since when we rotate them to the vertical position we obtain
crossings of the opposite type in the local tangles. In order to bear in mind this 
change of the local signs we put on the geometric picture the minuses on the even terms.  
 We shall use both ways of representation for extracting the
properties of rational knots and tangles.

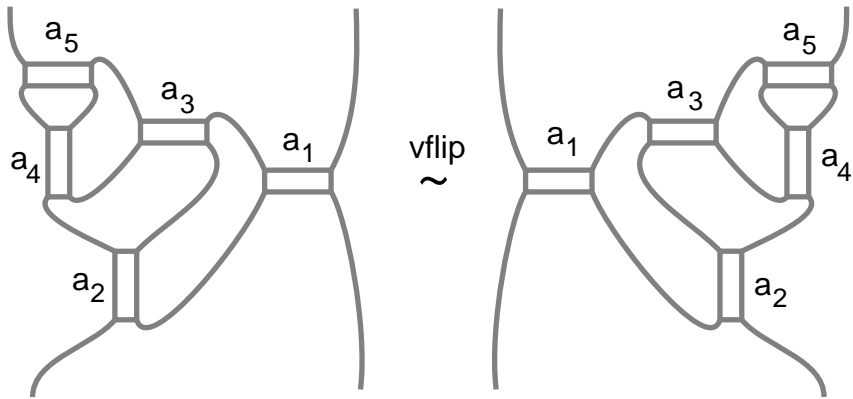
\bigbreak

$$ \picill4.7inby2.15in(R8) $$

\begin{center}
{Figure 8 - The standard representations } 
\end{center}
\vspace{3mm} 


\bigbreak

$$ \picill5.5inby1.9in(R9) $$

\begin{center}
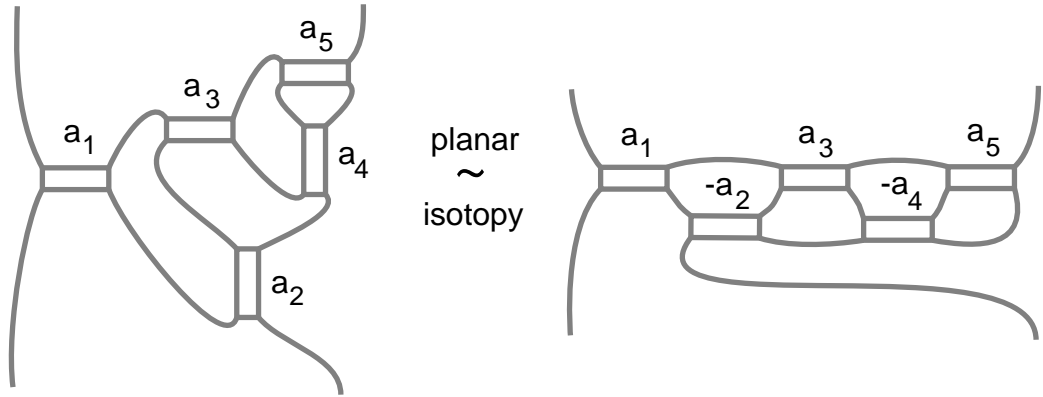
{Figure 9 - The standard and the $3$-strand-braid representation } 
\end{center}
\vspace{3mm} 


From the above one may associate to a
rational tangle diagram in standard form a vector  of integers $(a_1, a_2, \ldots, a_n), $
where the first entry denotes the place where the tangle starts unravelling and the last
entry where it  begins to twist.  For example the tangle of Figure 1 is associated  to
the vector
 $(2, -2, 3),$  while the  tangle of Figure 3 corresponds after a sequence of flypes to
the vector $(2, -4, -1, 3, 3).$ The vector associated to a
rational tangle diagram  is {\it unique} up to breaking the entry
 $a_n$ by a unit, i.e. $(a_1, a_2, \ldots, a_n) =  (a_1, a_2,
\ldots, a_n -1, 1),$ if $a_n >0,$ and $(a_1, a_2, \ldots, a_n) = (a_1, a_2, \ldots,
a_n +1, -1),$ if $a_n <0.$ This follows from the ambiguity of the very first crossing, see
Figure 10.  If  a rational tangle changes by an isotopy, the associated vector 
 might also change. 

\bigbreak

$$ \picill5inby.9in(R10) $$

\begin{center}
{Figure 10 - The ambiguity of the first crossing  } 
\end{center}
\vspace{3mm} 


\begin{rem}\label{odd}{\rm \   The same ambiguity implies that the number  $n$ in  the
above notation may be  assumed to be {\it odd}. We shall make this assumption for proving
Theorems~\ref{Schubert1} and~\ref{Schubert2}. } \end{rem}

 The next thing to observe is that a rational tangle in standard form can be described
 algebraically by a continued fraction built from the integer tangles  $[a_1],$
 $[a_2],$ $\ldots, [a_n]$  with all numerators equal to~$1,$ namely by an expression of
 the type:

 \[ [[a_1],[a_2],\ldots ,[a_n]] \, := \, [a_1]+ \frac{1}{[a_2]+\cdots +
 \frac{1}{[a_{n-1}]
 +\frac{1}{[a_n]}}} \]

 \noindent for $a_2, \ldots, a_n \in \ZZ - \{0\}$ and $n$ even or odd. We allow  $[a_1]$  to be the
tangle $[0].$  This expression follows inductively from the equation 

$$ T* \frac{1}{[n]} = \frac{1}{[n] +
 \frac{1}{T}}. $$  

 Then a rational tangle is
said to be in {\it continued fraction form}. 
 For example, Figure 1 illustrates the rational tangle $[[2], [-2], [3]],$  while the
tangles of  Figure 8 and 9 all depict
 the abstract rational tangle $[[a_1], [a_2], [a_3], [a_4], [a_5]].$

\bigbreak
 The tangle equation $ T* \frac{1}{[n]} = \frac{1}{[n] + \frac{1}{T}}$  implies
also that the two simple algebraic operations: {\it addition of $[+1]$ or $[-1]$} and {\it
inversion}  between rational tangles  generate the whole class of rational tangles.
For $T=[[a_1], [a_2],\ldots ,[a_n]]$ the following statements
 are
 now straightforward.

 \vspace{.15in}

 \noindent $\begin{array}{lrcl}

  1. & T + [\pm 1] &  = & [[a_1 \pm 1], [a_2],\ldots ,[a_n]],   \\
 [1.8mm]

 2. & \frac{1}{T} &  = & [[0], [a_1], [a_2],\ldots ,[a_n]],   \\  [1.8mm]

 3. & -T &  =  & [[-a_1], [-a_2],\ldots ,[-a_n]]. \\ [1.8mm]

 4. & T &  =  & [[a_1], [a_2], \ldots, [a_n -1], [1]], \ \ \mbox{if $a_n >0,$ and}  \\
[1.8mm]

   &  T &  =  & [[a_1], [a_2], \ldots, [a_n +1], [-1]],  \ \ \mbox{if $a_n <0.$}
 \end{array}$

\bigbreak

  A tangle is said to be {\it alternating} if the crossings 
alternate from under to over as we go along any component or arc of the weave.  Similarly, a knot is {\it
alternating} if it possesses an alternating diagram. We shall
see that rational tangles and rational knots are alternating. 
 Notice that, according to the checkerboard shading (see Figure 2 and the corresponding
discussion), the only way the weave alternates is if any two adjacent crossings are of the
same type,  and this propagates to the whole diagram. Thus, {\it a tangle or a knot
diagram with all crossings of the same type is alternating}, and this characterizes
alternating tangle and knot diagrams. It is important to note that flypes preserve the
alternating structure. Moreover, flypes are the only isotopy moves needed in the
statement of the  Tait Conjecture for alternating knots. An important property of
rational tangles is now the following.

\smallbreak
\noindent {\it A rational tangle diagram in standard form can be
always isotoped to an alternating one.
 } 
\smallbreak

\noindent  The process is inductive on the number of crossings and the basic
isotopy move is illustrated in Figure 11, see \cite{KL1} for details. We
point out that this isotopy applies to rational tangles in standard form where all the
crossings are on the right and on the bottom. We shall say that
a  rational tangle $T=[[a_1], [a_2],\ldots, [a_n]]$ is
 in  {\it  canonical form} if $T$ is alternating and $n$ is odd. From Remark 1 
we can always assume $n$ to be odd, so in order to bring a rational tangle to the
canonical form we just have to apply the isotopy moves described in Figure 11.  Note that
$T$ alternating implies that the $a_i$'s are all of the same sign.  

\bigbreak

 $$ \picill5.5inby1.35in(R11) $$

\begin{center}
{Figure 11 - Reducing to the alternating form  } 
\end{center}
\vspace{3mm} 


 \noindent  The alternating nature of the rational
 tangles will be  very useful to us in classifying rational knots and links.
 It turns out from the classification of alternating knots that {\it two
 alternating  tangles are isotopic if and only if they differ by a sequence
 of flypes.} (See \cite{SuT}, \cite{MT}. See also \cite{Sa}.) It is easy to see that the
closure of an alternating rational tangle is an alternating knot. Thus we have: 
\smallbreak
 \noindent { \it  Rational knots are alternating, since they possess a
 diagram that is the closure of an alternating rational tangle. }
 \bigbreak

\noindent {\bf Continued Fractions and the Classification of Rational Tangles.} 
From the above discussion it makes sense to assign to a rational tangle in standard
form, $T=[[a_1], [a_2],\ldots ,[a_n]],$ for 
$ a_1 \in \ZZ, \ a_2, \ldots, a_n \in \ZZ - \{0\}$ and $n$
 even or odd, the continued fraction 
 \[ F(T) = [a_1, a_2, \ldots, a_n] := a_1+ \frac{1}{a_2+\cdots + \frac{1}{a_{n-1}
 +\frac{1}{a_n}}},  \]

 \noindent if $T \neq [\infty],$ and $F([\infty]) := \infty = \frac{1}{0},$
 as a formal expression. This rational number or infinity shall be called {\it the
fraction  of $T$.} The fraction is a topological invariant of the tangle $T.$  We explain
briefly below how to see this. 

\bigbreak

 The subject of continued  fractions is of perennial interest
 to mathematicians. See for example 
  \cite{Kh},  \cite{O}, \cite{Ko}, \cite{W}. In this paper we shall only
 consider continued fractions of the above type, i.e. 
 with all numerators equal to 1. As in the
 case of rational tangles we allow the term $a_1$  to be zero. Clearly, the two simple
algebraic operations {\it addition  of $+1$ or $-1$} and {\it inversion}
 generate inductively the whole class of continued fractions starting from
 zero.  For any rational number $\frac{p}{q}$ the following statements are really
 straightforward.

 \vspace{.1in}

 \noindent \ $1.  \ \ \mbox{there are} \ a_1 \in \ZZ, \ a_2, \ldots, a_n
 \in
 \ZZ - \{0\} \
 \mbox{such that} \ \frac{p}{q} = [a_1, a_2, \ldots, a_n],  $

 \noindent $\begin{array}{lrcl}

  2.& \frac{p}{q} \pm 1 &  =  & [a_1 \pm 1, a_2,\ldots ,a_n],   \\
 [1.8mm]

  3. & \frac{q}{p} &  =  & [0, a_1, a_2, \ldots, a_n],   \\  [1.8mm]

 4. & -\frac{p}{q} &  =  & [-a_1, -a_2, \ldots, -a_n]. \\  [1.8mm]

 5. & \frac{p}{q} &  =  & [a_1, a_2, \ldots, a_n -1, 1], \ \ \mbox{if $a_n >0,$
       and}  \\ [1.8mm]

   &  \frac{p}{q} &  =  & [a_1, a_2, \ldots, a_n +1, -1],  \ \ \mbox{if $a_n <0.$}
 \end{array}$
 \vspace{.1in}

\noindent Property $1$ above is a consequence of Euclid's algorithm, see for example  
\cite{Kh}. Combining the above we obtain the following properties for the tangle fraction.

 \vspace{.1in}

 \noindent $\begin{array}{lrcl}

  1. & F(T + [\pm 1]) &  = & F(T) \pm 1,  \\
 [1.8mm]

  2. & F(\frac{1}{T}) &  = & \frac{1}{F(T)},   \\  [1.8mm]

 3. & F(-T) &  =  & -F(T). \\
 \end{array}$
 \vspace{.1in}

The last ingredient for the classification of rational tangles is the following fact
about continued  fractions. 

\smallbreak
\noindent { \it  Every  continued fraction $[a_1, a_2, \ldots, a_n]$ can
be
 transformed to a unique canonical form $ [\beta_1, \beta_2, \ldots, \beta_m],$ 
where
 all $\beta_i$'s are positive or all negative integers and $m$ is odd.
  }
\smallbreak

 \noindent One way to see this is to evaluate the continued fraction and then apply
Euclid's algorithm, keeping all remainders of the same sign. There is also an algorithm
that can be applied directly to the initial continued fraction to obtain its canonical
form. This algorithm works in parallel with the algorithm for the
canonical form of rational tangles, see \cite{KL1} for details. 
\bigbreak

From  the Tait
conjecture for alternating rational tangles, from the uniqueness of the canonical form of
continued fractions and from the above properties of the fraction we derive that the
fraction not only is an isotopy invariant of rational tangles but it also classifies
rational tangles. This is the Conway Theorem. See \cite{KL1} for details of
the proof. For the isotopy type of a rational tangle with fraction $\frac{p}{q}$ we
shall use the notation $[\frac{p}{q}].$ Finally,  it is easy to see the following useful
result about rational tangles:

\smallbreak
 { \it   Suppose that $T+ [n]$ is a rational tangle, then $T$ is a rational tangle.
   }

 \section{The Classification of Unoriented Rational Knots}

 In this section we shall prove Schubert's theorem for unoriented rational knots. 
It is convenient to  say that reduced fractions $p/q$ and $p'/q'$ are arithmetically
equivalent, written  $p/q \sim p'/q',$ if $p = p'$  and either $qq' \equiv 1 \, mod \, p$
\ or \  $q~\equiv~q'\,~mod \,~p.$ We shall call two rational tangles {\it arithmetically
equivalent} if their fractions are arithmetically equivalent.  In this language,
Schubert's theorem  states that two unoriented  rational tangles close to form isotopic
knots if and only if they are arithmetically equivalent. 

\smallbreak
 We  only need to consider numerator closures of rational tangles,
since the denominator closure of a tangle $T$ is simply the numerator closure of its
rotate $-\frac{1}{T}$. From the discussions in Section 2 a rational tangle may be
assumed to be in continued fraction form and by Remark
\ref{odd}, {\it the length} of a rational tangle may be assumed to be {\it odd.}  A rational knot
is said to be  {\it  in standard form, in continued fraction form, alternating or in
canonical form} if it is the numerator closure of a rational tangle that is 
in standard form, in continued fraction form, alternating or in canonical form
respectively.  By the alternating property of rational knots we may assume all rational
knot diagrams to be  {\it alternating}.  The diagrams and the isotopies of the rational
knots are meant to take place in the $2$-sphere and not in the plane. 

 \bigbreak

 \noindent {\bf Bottom twists.} \ The simplest instance of two rational tangles being 
non-isotopic but having isotopic  numerators is adding a number of twists at the
 bottom of a tangle, see Figure 12. Indeed, let $T$ be a rational tangle and let  $T*1/[n]$ be the
tangle obtained from $T$ by adding $n$ bottom twists, for any $n \in \ZZ.$ 
We have  $N(T*1/[n]) \sim N(T),$ but $F(T*1/[n]) = F(1/([n] + 1/T)) = 1/(n + 1/F(T));$
so, if 

\[  F(T) = p/q, \]

\noindent then

\[  F(T*1/[n]) = p/(np + q), \]

\noindent thus the two tangles are not isotopic. If we set $np + q =
q'$ we have $q\equiv
 q'\, mod \, p,$ just as Theorem \ref{Schubert1} predicts.  

\bigbreak

 $$ \picill3.4inby1.75in(R12) $$

\begin{center}
{Figure 12 -Twisting the Bottom of a Tangle } 
\end{center}
\vspace{3mm} 


Reducing all possible  bottom twists of a rational tangle yields a rational tangle
 with fraction $\frac{P}{Q}$ such that 

\[  |P|>|Q|.  \]  

To see this, suppose that we are dealing with
$\frac{P}{Q'}$ with $P < Q'$ and  both $P$ and $Q'$ positive (we leave it to the reader
to fill in the details for $Q'$ negative).  Then

$$ \frac{P}{Q'} = \frac{1}{\frac{Q'}{P}} = \frac{1}{n + \frac{Q}{P}} = \frac{1}{n +
\frac{1}{\frac{P}{Q}}},$$

\noindent where 
$$ Q'=nP +Q \equiv Q \, mod \, P, $$

\noindent for $n$ and $Q$  positive and $Q < P.$ So, by the Conway Theorem, the
rational tangle $[\frac{P}{Q'}]$ differs from the tangle
$[\frac{P}{Q}]$ by $n$ bottom twists, and so $N([\frac{P}{Q'}]) \sim  N([\frac{P}{Q}]).$
Figure 13 illustrates an example of this arithmetics. 
Note that a tangle with fraction $\frac{P}{Q}$ such that 
$ |P|>|Q| $ always ends  with a number of horizontal twists. So, if $T=
[[a_1],[a_2], \ldots, [a_n]]$ then $a_1 \neq 0.$ If $T$ is in
twist form then $T$ will not have any top or bottom twists. 
We shall say that a rational tangle whose fraction satisfies the above inequality is in
{\it reduced form}. 

\bigbreak

 $$ \picill5.5inby2.55in(R13) $$

\begin{center}
{Figure 13 - Reducing the Bottom Twists} 
\end{center}
\vspace{3mm} 


The proof of  Theorem \ref{Schubert1} now proceeds  in two stages.  First, (in 3.1) we look for all
possible places where we could cut a rational knot  $K$ open  to a rational tangle,
and we show that all cuts that open $K$ to other rational
tangles give tangles arithmetically equivalent to the tangle $T$.  
 Second, (in 3.2) given two isotopic reduced alternating rational knot diagrams, we have to check
that the rational tangles that they open to are arithmetically equivalent. By the
solution to the Tait Conjecture these isotopic knot diagrams will differ by a sequence of
flypes. So we analyze what happens when a flype is performed on $K.$

\subsection{The Cuts}  

Let  $K$ be a rational knot that is the numerator closure of a rational tangle  $T.$ 
 Wewill look for all `rational'  cuts on $K.$ 
 In our study of cuts we shall assume that $T$ is in reduced canonical form. The more general case
where $T$ is in reduced alternating twist form is completely analogous and we make a
remark at the end of the subsection. Moreover, the cut analysis in the case where
$a_1 = 0$ is also completely analogous for all cuts with appropriate adjustments. 
 There are three types of rational cuts.


\bigbreak

\noindent {\bf The standard cuts.} \  The tangle $T= [[a_1],[a_2], \ldots, [a_n]]$ is
said to arise as  {\it the standard cut} on $K=N(T).$ If we cut $K$ at another pair of
`vertical' points that are adjacent to the $i$th crossing of the elementary tangle $[a_1]$
(counting from the outside towards the inside of $T$) we obtain the alternating rational
tangle in twist form $T' =  [[a_1 - i],[a_2],
\ldots, [a_n]] + [i].$  Clearly, this tangle is isotopic to $T$ by a sequence of flypes that
send all the horizontal twists to the right of the tangle. See the right hand illustration of 
Figure 14 for $i=2.$ Thus, by the Conway Theorem, $T'$
will have the same fraction as $T.$  Any such cut on $K$ shall be called  a {\it standard cut} on
$K.$  

\bigbreak

$$ \picill4.8inby4.4in(R14) $$

\begin{center}
{Figure 14 - Standard Cuts  } 
\end{center}
\vspace{3mm} 


 \noindent {\bf The special cuts.} \ A key example of the arithmetic relationship of
the
 classification of rational knots is illustrated in Figure 15.
 The two tangles
 $T  = [-3] $ and $S = [1] + \frac{1}{[2]} $ are non-isotopic by the Conway
 Theorem, since $F(T) =  -3 = 3/-1, $ while $F(S) =1 +
 1/2 = 3/2.$ But they have isotopic
 numerators:  $N(T) \sim N(S),$ the
 left-handed trefoil. Now  $-1 \equiv 2 \, mod \, 3,$  confirming
 Theorem \ref{Schubert1}.

\bigbreak

$$ \picill3.45inby1.7in(R15) $$

\begin{center}
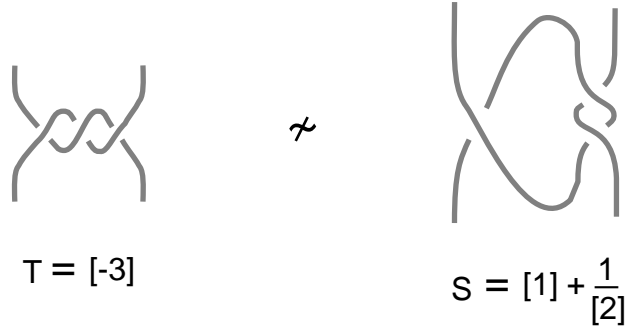
{Figure 15 - An Example of the Special Cut } 
\end{center}
\vspace{3mm} 


We now analyse the above example in general. Let $K = N(T),$ where $T=[[a_1],[a_2], \ldots, [a_n]].$
Since $T$ is assumed to be  in reduced  form, it follows that $a_1 \neq 0,$ so $T$ can be written in
the form
$T=[+1]+R$ \ or \ $T=[-1]+R,$ and the tangle $R$ is also rational. 

The indicated horizontal crossing $[+1]$ of the tangle $T=[+1]+R,$ which is the first
crossing of $[a_1]$ and the last created crossing of $T,$ 
 may also be seen as a vertical one. So, instead of cutting  the diagram $K$ open 
 at the two standard cutpoints to obtain the tangle $T,$  we  cut at the two other marked
`horizontal'
 points on the first crossing of the subtangle $[a_1]$ to obtain a new $2$-tangle
$T'$ (see Figure 16).  $T'$ is clearly rational, since $R$ is
rational. The tangle $T'$ is said to arise as  {\it the special cut} on $K.$ 
\smallbreak

\noindent We would like to identify this rational tangle
$T'.$ For this reason we first swing the upper arc of $K$ down to the bottom of the
diagram in order to free the region of the cutpoints. By our convention for the signs of crossings
in terms of the checkerboard shading, this forces all
crossings of $T$ to change sign from positive to negative and vice versa. We then rotate
$K$ by $90^0$ on its plane  (see right-hand illustration of Figure 16). This forces all
crossings of $T$ to change from horizontal to vertical and vice versa. In particular, 
the marked crossing $[+1],$ that was seen as a vertical one in $T,$ will now look as a
horizontal $[-1]$ in $T'.$ In fact, this will be the only last horizontal crossing of
$T',$ since all other crossings of
$[a_1]$ are now vertical.  So, if  $T=[[a_1],[a_2], \ldots, [a_n]]$ then  $R=[[a_1 - 1],
[a_2], \ldots ,[a_n]]$  and

\[ T' = [[-1], [1 - a_{1}], [-a_{2}],\ldots , [-a_{n}]].\] 

\noindent Note  that if the crossings of $K$ were all of negative type, thus  all the
$a_{i}$'s would be negative, the tangle $T'$ would be \  
$T' = [[+1], [-1 - a_{1}], [-a_{2}],\ldots , [-a_{n}]].$  In the example of Figure 15 if we took
$R= [-2], $ then $T= [-1] + R$ and $T' = S = [[+1], [+2]].$ 

\bigbreak

$$ \picill4.95inby4.75in(R16) $$

\begin{center}
{Figure 16 - Preparing for the Special Cut } 
\end{center}
\vspace{3mm}

The special cut is best illustrated in Figure 17. We consider the
rational knot diagram $K = N([+1]+R).$ (We analyze $N([-1]+R)$ in the same way.)  
As we see here, a
sequence of isotopies and  cutting $K$ open allow us to read the new tangle:  

\[ T' = [-1] - \frac{1}{R}. \]

\bigbreak

$$ \picill8inby3.1in(R17) $$

\begin{center}
{Figure 17 - The Tangle of the Special Cut } 
\end{center}
\vspace{3mm} 


From the above we have $N([+1]+R) \sim N([-1] -\frac{1}{R}).$  Let  now the fractions of
$T, R$ and $T'$ be $F(T) = P/Q, \, F(R) = p/q$ and $F(T') = P'/Q'$
respectively. Then

 $$F(T) = F([+1]+R) = 1+ p/q = (p+q)/q = P/Q,$$

\noindent  while 

$$F(T') = F([-1] - 1/R) = -1 - q/p = (p+q)/(-p) = P'/Q'.$$

\noindent The two fractions are different, thus the two rational tangles that give rise to
the same rational knot are not isotopic. We observe that $P=P'$ and 

$$ q \equiv -p \, mod  (p+q) \ \Longleftrightarrow \ Q  \equiv Q' \, mod \, P.$$

\noindent This arithmetic equivalence demonstrates another case for Theorem
\ref{Schubert1}. Notice that, although both the bottom twist and the special cut fall
into the same arithmetic equivalence, the arithmetic of the special cut is more subtle 
than the arithmetic of the bottom twist.

\bigbreak
 If we cut $K$ at the two lower horizontal points of the first crossing of 
$[a_1]$ we obtain the same rational tangle $T'.$ Also, if we cut at any other pair of
upper or lower horizontal adjacent points of the subtangle
$[a_1]$ we obtain a rational tangle in twist form isotopic to $T'.$ Such a cut shall be
called a {\it special cut}. See Figure 18 for an example.
 Finally, we may cut $K$ at any  pair of upper or lower horizontal adjacent points of the
subtangle $[a_n].$  We shall call this a {\it special palindrome cut.} We
will discuss this case after having analyzed the last type of a cut, the palindrome cut. 

\bigbreak

 $$ \picill5.5inby2.9in(R18) $$

\begin{center}
{ Figure 18 - A Special Cut} 
\end{center}
\vspace{3mm}

\noindent {\bf Note} \ We would like to point out that the horizontal-vertical ambiguity
of the last crossing of a rational tangle $T=[[a_1],\ldots ,[a_{n-1}], [a_n]],$  which
with the special cut on $K=N(T)$ gives rise to the tangle $[[\mp 1], [\pm 1 - a_1]
[-a_2],\ldots , [-a_n]]$   is very similar to the
horizontal-vertical ambiguity of the first crossing that does not change the tangle and
it gives rise to the tangle continued fraction $[[a_1],\ldots ,[a_{n-1}], [a_n \mp 1],
[\pm 1]].$

\begin{rem} \rm \ A special cut is nothing more than the addition of a bottom twist. Indeed,
 as Figure 19 illustrates, applying a positive bottom twist to the tangle $T'$ of the
special cut yields the tangle 
$S=([-1] -1/R)*[+1]$, and we find that if $F(R) = p/q$ then $F([+1] + R) = (p+q)/q$
while
$F(([-1] - 1/R)*[+1]) =  1/(1 + 1/(-1-q/p)) = (p+q)/q.$  Thus we see that the fractions of
$T=[+1]+R$ and $S=([-1] - 1/R)*[+1]$ are equal and by the Conway Theorem the tangle $S$ is 
isotopic to the original tangle $T$ of the standard cut. The isotopy move is nothing but
the transfer move of Figure 11. The isotopy is illustrated in Figure 19. Here we used
the Flipping Lemma. 
\end{rem}

\bigbreak

$$ \picill5.5inby1.5in(R19) $$

\begin{center}
{Figure 19 - Special Cuts and Bottom Twists  } 
\end{center}
\vspace{3mm}

 \noindent {\bf The palindrome cuts.} \  In Figure 20 we see that the tangles 

$$ T=[[2],[3],[4]] = [2] + \frac{1}{[3] + \frac{1}{[4]}} $$

\noindent and 
$$ S=[[4],[3],[2]] = [4] + \frac{1}{[3] + \frac{1}{[2]}} $$

\noindent both have the same numerator closure. This is another key example of the basic
relationship given in the classification of  rational knots. 

\smallbreak
In the general case if 
$ T=[[a_1],[a_2], \ldots, [a_n]], $ we shall call the tangle $ S=[[a_n],[a_{n-1}], \ldots,
[a_1]]$ {\it the  palindrome of $T$.} Clearly these tangles have the same numerator: 
$K = N(T) = N(S).$ Cutting open $K$ to yield $T$ is the standard cut, while
cutting to yield  $S$ shall be called  {\it the palindrome cut} on $K.$ 

\bigbreak

$$ \picill8inby2.5in(R20) $$

\begin{center}
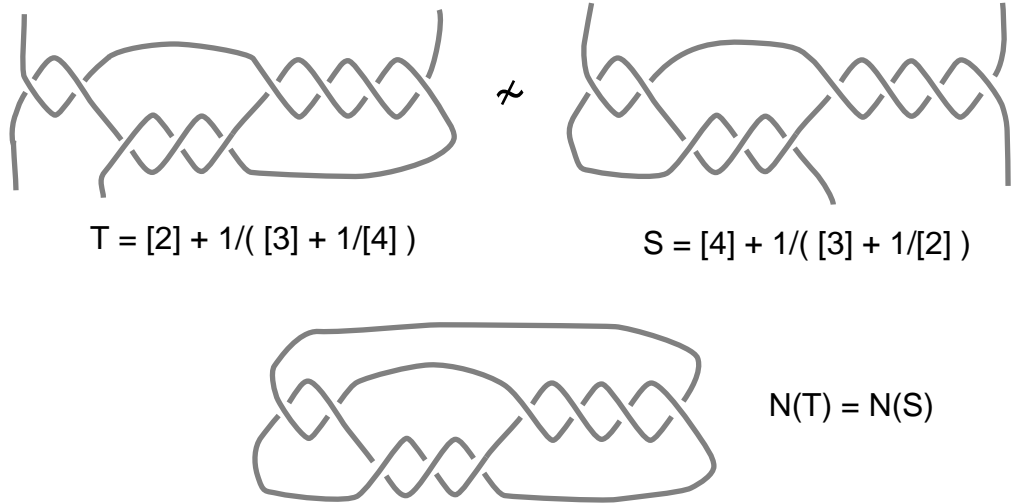
{Figure 20 - An Instance of the Palindrome Equivalence  } 
\end{center}
\vspace{3mm} 


The tangles in Figure 20 have corresponding fractions 
 $$ F(T)= 2 + \frac{1}{3 + \frac{1}{4}} = \frac{30}{13} \mbox{ \ \ and \ \  } F(S) = 4 +
\frac{1}{3 + \frac{1}{2}} = \frac{30}{7}.$$
 Note that \  $7\cdot 13\equiv 1 \, mod \, 30. $ This is the other instance of the
arithmetic behind the classification  of rational knots in Theorem \ref{Schubert1}. 
In order to check the arithmetic in the general case of the palindrome cut  we need to 
generalize this pattern to arbitrary continued fractions and their palindromes (obtained
by reversing the  order of the terms).  Then we have the following.

\begin{th}[Palindrome Theorem]\label{palindrome}{ \ Let $\{ a_{1},a_{2},\ldots ,a_{n} \}$
be a collection  of $n$   non-zero integers, and let $ \frac{P}{Q} = [a_{1},a_{2},\ldots
,a_{n}]$ and
$ \frac{P'}{Q'} =[a_{n},a_{n-1},\ldots ,a_{1}].$ Then 
$P=P'$ and  $QQ' \equiv  (-1)^{n+1}\, mod \, P.$
 } \end{th}

The Palindrome Theorem is a known result about continued fractions. For example see
\cite{Sie} or \cite{Kaw}, p. 25, Exercise 2.1.9. We shall give here our proof of this
statement. For this we will first present a way of evaluating continued fractions via
$2\times 2$ matrices (compare with \cite{Fr}, \cite{Ko}). This method of
evaluation is crucially important in our work in the rest of the paper. Let
$\frac{p}{q} = [a_{2},a_{3},\ldots ,a_{n}].$ Then we have:  

$$[a_{1},a_{2},\ldots ,a_{n}] = a_{1} + \frac{1}{\frac{p}{q}} = a_{1} +\frac{q}{p} =
\frac{a_1p+q}{p} = \frac{p'}{q'}.$$

\noindent Taking the convention that 
$[\left( \begin{array}{cc}
p\\
q\\
\end{array} \right)]
 := \frac{p}{q},$ with our usual conventions for formal fractions such as $\frac{1}{0},$ 
we can thus  write a corresponding matrix equation in the form

$$ \left( \begin{array}{cc}
a_{1} & 1 \\
1 & 0\\
\end{array} \right)
\cdot
\left( \begin{array}{cc}
p\\
q \\
\end{array} \right) 
=
\left( \begin{array}{cc}
a_{1}p+q\\
p\\
\end{array} \right) 
=
\left( \begin{array}{cc}
p'\\
q' \\
\end{array} \right). 
$$

\noindent We let

$$M(a_{i}) =
\left( \begin{array}{cc}
a_{i} & 1 \\
1 & 0\\
\end{array} \right).$$

\noindent  The matrices $M(a_i)$ are said to be the {\it generating matrices} for
continued  fractions,  as  we have:

\begin{lem}{\bf (Matrix interpretation for continued fractions)}{ \ For any sequence of non-zero
integers
 $\{  a_1,a_2,\ldots ,a_n \}$ 
the value of the corresponding continued fraction is given through the following  matrix
product

 $$[a_{1},a_{2},\ldots ,a_{n}] = [M(a_{1})M(a_{2})\cdot\cdot\cdot M(a_{n})\cdot v]$$

\noindent where 
$$v=
\left( \begin{array}{cc}
1\\
0\\
\end{array} \right).$$

 } \end{lem}

\noindent {\em Proof.}  We observe that 

\[ [M(a_{n}) \left( \begin{array}{cc}
1\\
0\\
\end{array} \right)] = 
[\left( \begin{array}{cc}
a_n\\
1\\
\end{array} \right)] = a_n = [a_n]  
\] 

and 

\[ [ M(a_{n-1}) \left( \begin{array}{cc}
a_n\\
1\\
\end{array} \right)] = [ \left( \begin{array}{cc}
a_{n-1} a_n + 1\\
a_n \\
\end{array} \right)] = [a_{n-1}, a_n]. \]

Now the lemma follows by induction. $\hfill \Box$

\bigbreak

\noindent {\em Proof of the Palindrome Theorem.} We wish to compare  $\frac{P}{Q} =
[a_{1},a_{2},\ldots ,a_{n}]$ and  $\frac{P'}{Q'} =[a_{n},a_{n-1},\ldots ,a_{1}].$ By 
Lemma 1 we can write
\vspace{1.8mm}

\noindent $\frac{P}{Q}=[M(a_{1})M(a_{2})\cdots M(a_{n})\cdot v]  \mbox{ \ \ and \ \ }
\frac{P'}{Q'} =[M(a_{n})M(a_{n-1})\cdots M(a_{1})\cdot v].$

\vspace{1.8mm}
\noindent Let 
$$M = M(a_{1})M(a_{2})\cdots M(a_{n})$$ 
and  
$$M' = M(a_{n})M(a_{n-1})\cdots M(a_{1}).$$ 
\vspace{1.8mm}

\noindent Then $ \frac{P}{Q}=[M\cdot v] $ and  $ \frac{P'}{Q'} =[M' \cdot v].$ We
observe that 
\vspace{1.8mm}

\noindent $M^T = (M(a_{1})M(a_{2})\cdots M(a_{n}))^T = (M(a_{n}))^T
(M(a_{n-1}))^T\cdots (M(a_{1}))^T $
\vspace{1.8mm}

$ \ =  M(a_{n})M(a_{n-1})\cdots M(a_{1}) = M',$  
\vspace{1.8mm}

\noindent since $M(a_{i})$ is symmetric, where 
$M^T$ is the transpose of $T.$ Thus 
$$ M' = M^T.$$ 

Let  
$$M =
\left( \begin{array}{cc}
X & Y \\
Z & U \\
\end{array} \right).$$
In order that the equations $ [M\cdot v] = \frac{P}{Q} $ \ and \ $ [M^T \cdot v] =
\frac{P'}{Q'} $ are satisfied it is necessary that  $X=P,$ $X=P',$ $Z=Q$ and $Y=Q'.$ 
That is, we should have:

$$M =
\left( \begin{array}{cc}
P & Q' \\
Q & U \\
\end{array} \right) \ \ \mbox{and} \ \ 
M' =
\left( \begin{array}{cc}
P & Q \\
Q' & U \\
\end{array} \right).$$

\noindent Furthermore, since  the determinant of $M(a_{i})$ is equal to $-1,$ we have
that 
$$det(M) = (-1)^n.$$
Thus 
$$PU - QQ' = (-1)^n,$$
 so that
 $$QQ' \equiv (-1)^{n+1}\, mod \, P,$$ 
and the proof of the Theorem is complete. $\hfill \Box$

\begin{rem}{\rm \ Note in the argument above that the entries of the matrix $M =
\left( \begin{array}{cc}
P & Q' \\
Q & U \\
\end{array} \right)$ of
a given continued fraction $[a_{1},a_{2},\ldots ,a_{n}] = \frac{P}{Q}$ involve also the 
evaluation of its palindrome  
$[a_{n},a_{n-1},\ldots ,a_{1}] = \frac{P}{Q'}.$  
 }
\end{rem}

Returning now to the analysis of the palindrome cut, we apply Theorem \ref{palindrome} in order to
evaluate the fraction of  palindrome rational tangles $T=[\frac{P}{Q}]$ and
$S=[\frac{P'}{Q'}].$ From the above analysis we have $P=P'.$ Also, by our assumption these
tangles  have continued fraction forms with odd length
$n,$ so we have the relation 

 $$QQ' \equiv 1\, mod \, P$$ 

\noindent and this is the second of the arithmetic relations of Theorem \ref{Schubert1}.
\smallbreak

If we cut  $K=N(T)$ at any other pair of `vertical' points of the subtangle
$[a_n]$ we obtain a rational tangle in twist form isotopic to the palindrome tangle $S.$  
 Any such cut shall be called a {\it palindrome cut}. 
\bigbreak

Having analysed the arithmetic of the palindrome cuts we can now return to the special
palindrome cuts on the subtangle $[a_n].$ These may be considered as special cuts on the
palindrome tangle $S.$ So, the fraction of the tangle of such a cut will satisfy the first
type of arithmetic relation of Theorem \ref{Schubert1} with the fraction of $S,$ namely
a relation of the type $q\equiv q'\, mod \, p,$  which, consequently, satisfies the second
type of arithmetic relation with the fraction of $T,$ namely a relation of the type 
$qq'\equiv 1\, mod \, p.$ In the end a  special palindrome cut will satisfy an arithmetic 
relation of the second type. This concludes the arithmetic study of the rational cuts.

\bigbreak

We now claim that the above listing of the three types of rational cuts is a complete
catalog of cuts that can  open the  link $K$ to a rational tangle: the standard cuts, the
special cuts and the palindrome cuts.  This is the crux of our  proof.

In Figure 21 we illustrate one example of a cut that is not rational. This is 
a possible cut made in the middle of the representative diagram $N(T)$. Here we see 
that if $T'$ is the tangle obtained from this cut, so that $N(T') = K$, then $D(T')$ is a
connected sum of two non-trivial knots. Hence the denominator $K'=D(T')$ is not prime.
Since we know that both the numerator and the denominator of a rational tangle are prime
(see 
\cite{BZ}, p. 91 or \cite{Li}, Chapter 4, pp. 32--40), it follows that $T'$ is not
a rational tangle. Clearly the above argument is generic. It is not  hard to see by
enumeration that all possible cuts with the exception of the ones we have described will
not give rise to rational tangles. We omit the enumeration of these cases.  

\smallbreak
This completes the proof that all of the rational tangles that close to a given standard rational
knot diagram are arithmetically equivalent. 

\bigbreak

$$ \picill8inby2.3in(R21) $$

\begin{center}
{ Figure 21 - A Non-rational Cut} 
\end{center}


In Figure 22 we illustrate on a representative rational knot in $3$-strand-braid form all
the cuts that exhibit that knot as a closure of a rational tangle. Each pair of points is
marked with the same number. 

\bigbreak

$$ \picill8inby2.85in(R22) $$

\begin{center}
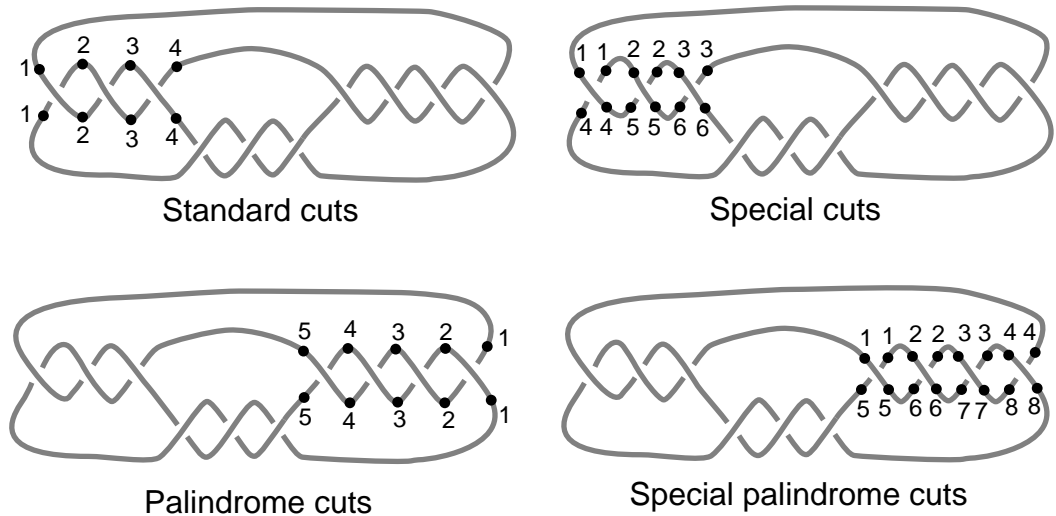
{ Figure 22 - Standard, Special, Palindrome and Special Palindrome Cuts } 
\end{center}
\vspace{3mm}


\begin{rem}{\rm \ It follows from the above analysis that if $T$ is a rational tangle in 
twist form, which is isotopic to the standard form $[[a_1],[a_2], \ldots, [a_n]],$ then
all arithmetically equivalent rational tangles can arise by any cut of the above types
either on the crossings that add up to the subtangle   $[a_1]$ or on the
crossings of the subtangle $[a_n].$ }
\end{rem}

\subsection{The flypes}

Diagrams for knots and links are represented on the surface of the two-sphere, $S^2,$ and
then notationally on a plane for purposes of illustration. 

Let $K=N(T)$ be  a rational link diagram with $T$ a rational tangle in twist form. 
By an appropriate sequence of flypes (recall Definition 1)  we may assume, without loss of
generality, that $T$ is alternating and in continued fraction form, i.e. $T$ is of the
form $T= [[a_1],[a_2], \ldots, [a_n]]$ with  the $a_i$'s all positive or all
negative. From the ambiguity of the first crossing of a rational tangle we may  assume
that $n$ is odd.  Moreover, from the analysis of the bottom twists in the previous subsection   we
may  assume that $T$ is in reduced form. Then the numerator $K=N(T)$  will be a reduced 
alternating knot diagram. This follows from the primality of $K.$ 

Let $K$ and $K'$ be two
 isotopic, reduced, alternating rational knot diagrams. By the Tait Conjecture they
will differ by a finite sequence of flypes. In considering  how rational knots can be cut
open to produce  rational  tangles, we will examine how the cuts are affected by
flyping. We analyse all
possible flypes to prove that it is sufficient to consider the cuts on a
single alternating reduced diagram for a given rational knot $K.$  Hence the
proof will be complete at that point. We need first two definitions and an observation
about flypes.

\begin{defn}{\rm \ We shall call {\it  region of a flype} the part of the knot diagram
that contains precisely the subtangle and the crossing that participate in the flype. The
region of a flype can be enclosed by a simple closed curve on the plane that intersects
the tangle in four points. }
\end{defn}

\begin{defn}{\rm \ A {\it pancake flip} of a knot diagram in the plane is an isotopy move
 that rotates the diagram by  $180^0$ in space around a horizontal or vertical
axis on its plane and then it replaces it on the plane. Note that any knot diagram in
$S^2$ can be regarded as a knot diagram in a plane. }
\end{defn}

In fact, the pancake flip is actually obtained by flypes so long as we
allow as background moves isotopies of the diagram in $S^2.$
To see this, note as in Figure 23 that we can assume without loss of
generality that the diagram in question is of the form 
$N([\pm 1] + R)$ for some tangle $R$ not necessarily rational. 
(Isolate one crossing at the `outer edge' of the diagram in the plane and
decompose the diagram into this crossing and a complementary tangle, as
shown in Figure 23.)  In order to place the
diagram in this form we only need to use isotopies of the diagram in the plane. 

\bigbreak

$$\vbox{\picill5.5inby1.8in(R23)  }$$

\begin{center}
{ Figure 23 - Decomposing into $N([\pm 1] + R)$} 
\end{center}
\vspace{3mm}

  Note now, as in Figure 24, that the pancake flip applied to
$N([\pm 1] + R)$  yields a diagram that can be obtained by a combination of a
planar isotopy, $S^2$-isotopies and a flype. (By an {\it $S^2$-isotopy} we mean the
sliding of an arc around the back of the sphere.) This is valid for $R$ any
$2$-tangle. We will use this remark in our study of  rational knots and links.
 
\bigbreak

$$\vbox{\picill5.5inby2.55in(R24)  }$$

\begin{center}
{ Figure 24 - Pancake Flip} 
\end{center}
\vspace{3mm}

We continue with a general remark about the form of a flype in any knot or link in $S^2.$ 
View Figure 25.  First look at parts A and B of this figure. Diagram A
is shown as a composition of a crossing and two tangles $P$ and $Q.$ 
 Part B is obtained from a flype of part A, where the flype
is applied to the crossing in conjunction with the tangle $P.$ 
This is the general pattern of the application of a flype. The flype is
applied to a composition of a crossing with a tangle, while the rest of 
the diagram can be regarded as contained within a second tangle that is
left fixed under the flyping.  

Now look at diagrams C and D. Diagram D is obtained by a flype using $Q$ and a crossing on
diagram C.  But diagram C is isotopic by a planar isotopy to diagram A, and diagrams B
and D are related by a pancake flip (combined with an isotopy that
swings two arcs around $S^2$). Thus we see that: 

\bigbreak
\noindent {\em Up to a pancake flip one can choose to
keep either of the tangles $P$ or $Q$ fixed in performing a flype}.

\bigbreak

$$\vbox{\picill5.5inby2.95in(R25)  }$$

\begin{center}
{ Figure 25 - The Complementary Flype} 
\end{center}
\vspace{3mm}

Let now $K=N(T)$ and $K'=N(T')$ be two reduced alternating rational knot diagrams that
differ by a flype. The rational tangles $T$ and $T'$ are in reduced alternating twist
form and without loss of generality $T$ may be assumed to be in continued fraction form. 
Then, recall from Section 2 that the region of the flype on $K$ can either include a rational
truncation of $T$ or some crossings of a subtangle $[a_i],$  see Figure
26. In the first case the two subtangles into which $K$ decomposes are both rational
and each will be called the {\it complementary tangle} of the other. In the second
case the flype has really trivial effect and the complementary tangle is not rational,
unless
$i = 1
\mbox{ or } n.$ 

\bigbreak

$$\vbox{\picill8inby1.5in(R26)  }$$

\begin{center}
{ Figure 26 -  Flypes of Rational Knots} 
\end{center}
\vspace{3mm}

For the cutpoints of $T$ on $K=N(T)$ there are three possibilities: 

\bigbreak
1.   \ they are outside the region of the flype, 
\smallbreak
2.  \ they are inside the flyped subtangle, 
\smallbreak
3.  \ they are inside the region of the flype and outside the flyped subtangle. 
\bigbreak

 If the cutpoints are outside the region of the flype, then the flype is
taking place inside the tangle $T$ and so there is nothing to check, since the new
tangle is isotopic and thus arithmetically equivalent to $T.$ 

\smallbreak
We concentrate now on the first case of the region of a flype. If the cutpoints
are inside the flyped subtangle then, by Figure 25, this flype can be seen as a
flype of the complementary tangle followed by a pancake flip. The region of the flype of
the complementary tangle does not contain the cut points, so it is a rational flype that
isotopes the tangle to itself.  The pancake flip also does not affect the arithmetic,
because its effect on the level of the tangle $T$ is simply a horizontal or a vertical
flip, and we know that a flipped rational tangle is isotopic to itself. 

If the region of the flype encircles a number of crossings of some $[a_i]$ then the
cutpoints cannot lie in the region, unless $i = 1 \mbox{ or } n.$  If the cutpoints do
not lie in the region of the flype, there is nothing to check. If  they do, then 
the complementary tangle is isotopic to $T,$ and the pancake flip produces an isotopic
tangle.

\smallbreak
Finally, if the cutpoints are inside the region of the flype and outside the flyped
subtangle, i.e. they are near the crossing of the flype, then there are three cases to
check. These are illustrated in Figure 27.  

\bigbreak

$$ \picill4.75inby2.3in(R27) $$

\begin{center}
{ Figure 27 - Flype and Cut Interaction} 
\end{center}

  In each of these cases the flype is illustrated with respect to a crossing and a tangle 
$R$ that  is a subtangle of the link $K=N(T).$ Cases (i) and (ii) are taken care of by the
trick of the complementary flype. Namely, as in Figure 25, we transfer the crossing
 of the flype around $S^2.$ Using this crossing we do a tangle flype of the complementary
tangle, then we do a horizontal pancake flip and finally an $S^2$-isotopy, to end up with
the right-hand sides of Figure 27. 

In case (iii) we note that after the flype
the position of the cut points is outside the region of a flyping move that can be
performed on the diagram $K'$ to return to the original diagram $K,$ see Figure 28.  
This means that after performing the return flype  the tangle $T'$ is isotopic to the
tangle $T''.$  One can now observe that if the original cut 
produces  a rational tangle, then the cut after the returned flype also  produces a
rational tangle, and this is arithmetically equivalent to the tangle $T.$ 
 More precisely,  the tangle $T''$ is the result of a special cut on $N(T).$

\bigbreak

$$ \picill6inby.85in(R28) $$

\begin{center}
{ Figure 28 - Flype and Special Cut} 
\end{center}

With the above argument we conclude the proof of the main direction of Theorem
\ref{Schubert1}. From our analysis it follows that: 
\bigbreak
\noindent {\it If $K=N(T)$ is a rational knot
diagram with
$T$ a rational tangle then, up to bottom twists, any other  rational tangle that closes to
this knot is available as a cut on the given diagram. }
\bigbreak

 We will now show the converse. 
\smallbreak
\noindent
$``\Longleftarrow "$ \ We wish to show that if two rational tangles are arithmetically 
equivalent, then their numerators are isotopic  knots. Let $T_1, T_2$ be rational 
tangles with
$F(T_1) = \frac{p}{q}$ and 
$F(T_2) = \frac{p}{q'},$ with $|p|>|q|$ and $|p|>|q'|,$  and assume first $qq'\equiv 1 \,
mod \, p.$  If $\frac{p}{q} = [a_1, a_2, \ldots, a_n],$ with $n$ odd, and $\frac{p}{q''} = [a_n,
a_{n-1}, \ldots, a_1]$ is the corresponding palindrome  continued fraction, then it follows from 
the Palindrome Theorem that  $qq''\equiv 1 \, mod \, p.$ Furthermore, it follows  by induction
that in a product of the form
$$M(a_1)M(a_2) \cdots M(a_n) = \left( \begin{array}{cc}
p & q'' \\
q & u \\
\end{array} \right)$$ 
we have that $p>q$ and $p>q''$, $q \ge u$ and $q'' \ge u$ whenever $a_1, a_2, \cdots, a_n$ are positive
integers. (With the exception in the case of $M(1)$ where the first two inequalities are replaced by
equalities.) Induction step involves multiplying a matrix in the above form by one more matrix $M(a)$, and
observing that the inequalities persist in the product matrix. 

Hence, in our discussion we can conclude that $|p|>|q''|.$
Since $|p|>|q'|$ and $|p|>|q''|,$ 
it follows that $q'= q'',$ since they are both reduced residue solutions of a $mod \, p$ equation 
with  a unique solution. Hence  
$[a_n,a_{n-1},\ldots, a_1] = \frac{p}{q'},$  
and, by the uniqueness of the canonical form for rational tangles, $T_2$
has to be: 

$$T_2 = [[a_n], [a_{n-1}], \cdots ,[a_1]].$$ 
For these tangles we know that $N(T_1) = N(T_2).$ 
Let now $T_3$ be another rational tangle with fraction 

$$\frac{p}{q' + k \, p} = \frac{1}{\frac{q'}{p} + k}.$$ 
By the Conway Theorem, this is the fraction of the rational tangle  
$$\frac{1}{\frac{1}{T_2} + [k]} =  T_2 * \frac{1}{[k]}.$$ 
 Hence we have (recall the analysis of the
bottom twists):
$$N(\frac{1}{\frac{1}{T_2} + [k]}) \sim N(T_2).$$

Finally, let $F(S_1) = \frac{p}{q}$ and  $F(S_2) = \frac{p}{q + k \, p}.$ Then 
$$\frac{p}{q + k \, p} = \frac{1}{\frac{q}{p} + k},$$ which is the fraction of the
rational tangle  
$$\frac{1}{\frac{1}{S_1} + [k]} = S_1 * \frac{1}{[k]}.$$ Thus 
$$N(S_1) \sim N(S_2).$$  The proof of Theorem
\ref{Schubert1} is now finished. $\hfill \Box $

\bigbreak
We close the section with two remarks.

\begin{rem}{\rm \  In the above discussion about flypes we 
used the fact that the tangles and flyping tangles involved were rational. One
can consider the question of {\em arbitrary alternating    tangles $T$ that close to form
links isotopic to a given  alternating diagram $K$.} Our analysis of cuts occuring before
and after a flype goes through   to show that {\em for every alternating tangle
$T,$ that  closes to a diagram isotopic to a given alternating diagram $K$, there is a cut 
on the diagram $K$ that produces a tangle that is arithmetically equivalent to $T$.} Thus
it makes sense to consider the collection of tangles that close to an arbitrary
alternating link up to this arithmetic equivalence. In the general case of
alternating links this shows that on a given diagram of the alternating link we can
consider all cuts that produce alternating tangles and thereby obtain all such tangles,
up to a certain arithmetical equivalence, that close to links isotopic to
$K.$ Even for rational links there can be more than one equivalence class of such
tangles.  For example, $N(1/[3] + 1/[3]) = N([-6])$ and
$F(1/[3] + 1/[3]) = 2/3$ while $F([-6]) = -6.$ Since these fractions have different 
numerators their tangles (one of which is not rational) lie in different equivalence
classes.  These remarks lead us to consider the set of
arithmetical equivalence classes of  altenating tangles that close to a given 
alternating link and to search for an analogue of Schubert's Theorem in this general
setting. 
 }
\end{rem}

\begin{rem}{\rm \  DNA supercoils, replicates and recombines with the help of certain
 enzymes. {\it Site-specific
 recombination} is one of the ways nature alters the genetic code of an
 organism, either by moving
 a block of DNA to another position on the molecule or by integrating a
 block of alien DNA into a
 host genome.  In \cite{CSS} it was
 made  possible for the first time to see knotted DNA in an electron micrograph with
 sufficient resolution to actually
 identify the topological type of these knots and links. It was possible to
 design an experiment involving successive DNA recombinations and to examine
 the topology of the products.  In
 \cite{CSS} the knotted DNA  produced by such successive recombinations was
 consistent with the
 hypothesis that all recombinations were of the type of a positive half
 twist as in $[+1].$ Then
 D.W. Sumners  and C. Ernst \cite{ES2} proposed a {\em tangle model for
 successive DNA
 recombinations} and showed, in the case of the experiments in question,
 that there  was no other
 topological possibility for the recombination mechanism than the positive
 half twist $[+1].$ Their work depends essentially on the classification theorem for
rational knots. This
 constitutes a unique use of topological mathematics as a theoretical
 underpinning for a problem in molecular biology.
 }
\end{rem}

\section{Rational Knots and Their Mirror Images}

 In this section we give an application of Theorem \ref{Schubert1}.  
 An unoriented knot or link $K$ is said to be {\em achiral} if it is
 topologically equivalent to its
 mirror image $-K$. If a link is not equivalent to its mirror image then it
 is said be {\em chiral}.
 One then can speak of the  {\em chirality} of a given knot or link,
 meaning whether it is chiral or
 achiral. Chirality plays an important role in the applications of knot
 theory to chemistry and
 molecular biology. In \cite{ES1} the authors find an explicit formula for the number of
achiral rational knots among all rational knots with $n$ crossings. 
It is interesting to use the classification of rational knots and links to
 determine their chirality. Indeed, we have the following well-known result
(for example see  \cite{Sie} and \cite{Kaw}, p. 24, Exercise 2.1.4. Compare also with
\cite{Sch2}):

 \begin{th}\label{mirror}{\ Let $K=N(T)$ be an unoriented rational knot or link,
 presented as the numerator of a
 rational  tangle $T.$  Suppose that  $F(T) = p/q$ with $p$ and $q$
 relatively prime. Then $K$  is achiral if and only if
  $q^{2}\equiv -1 \, mod \, p.$  It follows that the tangle $T$ has to be of the form
 $[[a_1], [a_2], \ldots, [a_k], [a_k], \ldots, [a_2], [a_1]]$ for any
 integers $a_1, \ldots, a_k.$ } 
 \end{th}

\noindent  Note that in this description we are using a representation of the tangle with
an even number of terms. The leftmost  twists $[a_1]$ are horizontal, thus $|p|>|q|.$
The rightmost starting twists are then vertical. 
\smallbreak 

 \noindent {\em Proof.}   With $-T$ the mirror image of the tangle $T$, we
 have that $-K = N(-T)$ and $F(-T) = p/(-q).$ If $K$ is isotopic to $-K,$ it follows from
the  classification theorem for rational
 knots that either $q(-q) \equiv 1 \, mod \, p$ or $q \equiv -q \, mod \, p.$ Without
loss of generality we can assume that $0< q < p.$ Hence $2q$ is not divisible by $p$ and
therefore it is not the case that  $q \equiv -q \, mod \, p.$ Hence $q^{2} \equiv -1 \,
mod \, p.$
 \smallbreak

 \noindent Conversely, if $q^{2} \equiv -1 \, mod \, p,$ then it follows from the
Palindrome Theorem that {\it the
 continued fraction expansion of $p/q$ has to be palindromic with an even number of terms.}
To see this, let $p/q = [c_1, \cdots , c_n]$ with n even, and let $p'/q' = [c_n, \cdots , c_1].$
The Palindrome theorem tells us that $p'= p$ and that $qq' \equiv -1 \, mod \, p.$ Thus we have
that both $q$ and $q'$ satisfy the equation $qx \equiv -1 \, mod \, p$ and both $q$ and $q'$ are between
$1$ and $p-1.$ Since this equation has a unique solution in this range, we conclude that $q=q'.$
It follows at once that the continued fraction sequence for $p/q$ is symmetric.

It is then easy
 to see that the
 corresponding rational knot or link $K = N(T)$  is equivalent to  its mirror
image. One
 rotates $K$ by $180^0$ in the plane and swings an arc, as Figure 29 illustrates. 
The point is that the crossings of the second row of the tangle $T,$  that are
seemingly crossings  of opposite type than the crossings of the upper row, become
after the turn crossings of the  upper row, and so the types of crossings are
switched. 
 This completes the proof.
 $\hfill \Box$


\bigbreak

 $$ \picill4.5inby2in(R29) $$

\begin{center}
{ Figure 29 - An Achiral Rational Link } 
\end{center}
\vspace{3mm}


 \section{On connectivity}

 We shall now introduce  the notion of {\it connectivity} and we shall
 relate it to  the fraction of unoriented rational tangles. 
 We shall say that an unoriented rational tangle has {\it connectivity type $[0]$}
 if the  NW end arc is connected to the  NE end arc and the  SW  end
 arc is connected to the SE end arc. These are the same connections as in the
 tangle $[0]$. Similarly, we say that 
 the tangle has {\it connectivity type $[\infty]$} or {\it $[1]$ } 
 if the end arc connections
 are the same as in the tangles  $[\infty]$  and $[+1]$ (or equivalently $[-1]$)
respectively. 
 The  basic connectivity patterns of rational tangles are exemplified by the tangles
 $[0]$, $[\infty]$ and $[+1]$.  We can represent them iconically by

 $$[0] = \mbox{\large $\asymp$}$$
 $$[\infty] = ><$$
 $$[1] =\mbox{\large $\chi$}$$

\noindent For connectivity we are only concerned with the connection patterns of the four 
end arcs. Thus $[n]$ has connectivity $\mbox{\large $\chi$}$ whenever $n$ is odd, and
connectivity $\mbox{\large $\asymp$}$ whenever
$n$ is even. Note that connectivity type $[0]$ yields two-component rational links, 
whilst  type  $[1]$ or $[\infty]$ yields one-component rational links. Also, adding a
bottom  twist to a rational tangle of  connectivity type $[0]$ will not change the
connectivity type of the tangle, while  adding a bottom twist to a rational tangle
of  connectivity type $[\infty]$ will switch the connectivity type to $[1]$ and  
vice versa. 
\smallbreak

\noindent We need to keep an accounting of the  connectivity
 of rational  tangles in relation to the parity of the  numerators and
 denominators of their fractions. 
 We adopt the following notation:  $e$ stands for {\it even} and $o$ 
 for {\it odd}. The {\it parity of a fraction} 
 $p/q$ is defined to be the ratio of the parities ($e$ or $o$) of its numerator 
 and denominator $p$ and $q$. Thus the fraction $2/3$ is of parity $e/o.$ 
 The tangle $[0]$ has fraction $0 = 0/1,$ thus parity $e/o.$ The tangle 
 $[\infty]$ has formal fraction $\infty = 1/0,$ thus parity $o/e.$  
 The tangle $[+1]$ has fraction $1 = 1/1,$ thus parity $o/o,$ and the tangle $[-1]$ has
fraction $-1 = -1/1,$ thus parity $o/o.$ We then have the
 following result.

 \begin{th}\label{connectivity} \ A rational tangle $T$ has  connectivity type 
 $\mbox{\large $\asymp$}$ if and only if its fraction has parity $e/o$.
  $T$ has  connectivity type $><$ if  and only if its fraction has parity
 $o/e$. Finally, $T$ has  connectivity type $\mbox{\large $\chi$}$ if and only if
 its fraction has parity $o/o$. 
   \end{th}

\noindent {\em Proof.} Since $F([0]) = 0/1$, $F([\pm 1]) = \pm 1/1$ and
$F([\infty]) = 1/0$, the theorem is true for these elementary
tangles. It remains to show by induction that it is true for any rational tangle
$T$. Note how connectivity type behaves under the addition
and product of tangles:  

$$\mbox{\large $\chi$} + \mbox{\large $\chi$} = \mbox{\large $\asymp$}$$
$$\mbox{\large $\chi$} + \mbox{\large $\asymp$} =\mbox{\large $\chi$} $$
$$\mbox{\large $\asymp$} + \mbox{\large $\asymp$} = \mbox{\large
$\asymp$}$$
$$\mbox{\large $\chi$} + >< = ><$$
$$\mbox{\large $\asymp$} + >< = ><$$
$$>< + >< = > <> < = \delta ><$$
\smallbreak

$$\mbox{\large $\chi$} * \mbox{\large $\chi$} = >< $$
$$\mbox{\large $\chi$} * \mbox{\large $\asymp$} = \mbox{\large $\asymp$}$$
$$\mbox{\large $\asymp$} * \mbox{\large $\asymp$} = \delta \, \mbox{\large
$\asymp$}$$
$$\mbox{\large $\chi$} * >< = \chi$$
$$\mbox{\large $\asymp$} * >< = \mbox{\large $\asymp$} $$
$$>< * >< = > <$$
\smallbreak

\noindent The symbol $\delta$ stands for the value of a loop formed. Now any
rational tangle can be built from $[0]$ or
$[\infty]$ by successive addition or multiplication with $[\pm 1].$ Thus, from the 
point of view of connectivity, it suffices to show that $[T] + [\pm 1]$ and
$[T]*[\pm 1]$ satisfy the theorem whenever $[T]$ satisfies the theorem.
This is checked by comparing the connectivity identities above with the
parity of the fractions. For example, in the case 
 $$ \mbox{\large $\chi$} + \mbox{\large $\chi$} =
\mbox{\large
$\asymp$} \ \ \mbox{we have} \ \ o/o + o/o = e/o$$ 
\noindent  exactly in accordance with the
connectivity identity. The other cases correspond as well, and this 
proves the theorem by induction.
$\hfill \Box $

\begin{cor}{ \ For a rational tangle $T$ the link  $N(T)$ has two components if
 and only if $T$ has fraction $F(T)$ of parity $e/o.$
 } \end{cor} 

\noindent {\em Proof.} By the Theorem we have $F(T)$ has parity $e/o$ if
and only it $T$ has connectivity type 
$\mbox{\large $\asymp$}.$ It follows at once that $N(T)$ has two
components.
$\hfill \Box$
\bigbreak

 Another useful fact is that the components of a rational link
are individually unknotted embeddings in three dimensional space. 
 To see this, observe that upon removing one strand of a rational
tangle, the other strand is an unknotted arc.

 \section{The Oriented Case}

Oriented rational knots and links are numerator (and denominator) closures of oriented
rational tangles. Rational tangles are oriented by choosing an orientation
for each strand of the tangle. Two oriented rational tangles are {\it isotopic} if
they are isotopic as unoriented tangles via an isotopy that carries the orientation of
one tangle to the orientation of the other.  Since the end arcs of a tangle are fixed
during a tangle isotopy, this means that  isotopic tangles must have identical
orientations at their end arcs. Thus,   {\it two oriented tangles are isotopic if and
only if they are isotopic as unoriented tangles and  they have identical orientations at
their end arcs.}  It follows that a given unoriented rational tangle can always yield
non-isotopic oriented rational tangles, for different choices of orientation  of one or
both strands. 
\smallbreak
 In  order to compare oriented rational knots via rational tangles we are only
interested in orientations that  yield  consistently oriented knots  upon taking the
numerator closure. This means that the two top end arcs 
have to be oriented one inward and the other outward. Same for the two bottom  end
arcs. 

\smallbreak
Reversing the orientation of one strand of an oriented rational tangle that gives rise
to  a two-component link will usually yield non-isotopic oriented rational
links. Figure 30  illustrates an example of non-isotopic oriented rational links, which
are isotopic as unoriented links. But reversing a single strand may also yield isotopic oriented rational
links. This will be the subject of the next section. 

\bigbreak

$$ \picill8inby1.15in(R30) $$

\begin{center}
{ Figure 30 - Non-isotopic Oriented Rational Links } 
\end{center}
\vspace{3mm}


An oriented knot or link is said to be {\it invertible} if it is oriented isotopic to
its inverse, i.e. the link obtained from it by reversing the orientation of each
component. We can obtain the inverse of a rational link  by reversing the orientation of
both strands of the oriented rational tangle of which it is the numerator. 
 It is easy to see that any rational knot or link is invertible. See the example 
  on the right-hand side of Figure 31.

\begin{lem} \label{invertible} \ Rational knots and links are invertible.  
  \end{lem} 

\noindent {\em Proof.} Let $K = N(T)$ be an oriented rational knot or link with $T$ 
an oriented rational tangle.  We
do a  vertical $180^0$-rotation in space, as the left-hand side of Figure 31
demonstrates. This rotation is a vertical flip for the rational tangle $T.$ Let $T'$
denote the result of the vertical flip of the tangle $T.$  The resulting oriented
knot $K'=N(T')$ is  oriented isotopic to $K,$ while the orientation of $T'$ is
the opposite of that of $T$ on both strands, and thus on all end arcs. But as we have
already noted $T$ is  isotopic to its vertical flip as unoriented tangles, thus they will
have the same fraction. It follows that 
$T'$ can be isotoped to $T$ through an (unoriented) isotopy  that  leaves the external
strands fixed. Therefore, the result of the vertical $180^0$-rotation is the tangle $T$
but with all orientations reversed. Thus  $K'$ is the
inverse of $K,$ and from the above $K$ is oriented isotopic to its inverse.   $\hfill
\Box$
 
\bigbreak

 Using this observation we conclude that, as far as the
study of  oriented rational knots is concerned, {\it all oriented rational tangles  may
be assumed to have the same orientation for their two upper end arcs.} Indeed, if
the orientations of the two upper end arcs are opposite of the fixed ones we do a
vertical flip to obtain the orientation pattern that agrees with our convention. 
 We fix this orientation to be downward for the NW end arc and  upward for the  NE end
arc, as in the examples of Figure 30 and as illustrated in Figure 32. 
\bigbreak

$$ \picill4.75inby1.2in(R31) $$

\begin{center}
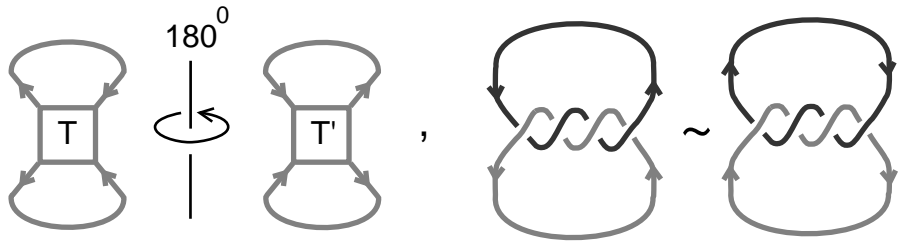
{ Figure 31 - Isotopic Oriented Rational Knots and Links } 
\end{center}
\vspace{3mm}


 \noindent  Thus we may reduce our analysis to two  basic
 types of orientation for the four end arcs of a rational tangle. 
 We shall call an oriented rational tangle {\it of type I} if the  SW arc is oriented
  downward and the   SE arc is oriented upward, and {\it of type II} if the SW arc
is oriented upward and the  SE arc is oriented downward, see Figure 32. From the
 above remarks  any tangle is of type
 I or type II. Two tangles are said to be {\it compatible} it they are both
 of type I or both of  type II and {\it incompatible} if they are of different types.
\smallbreak
\noindent {\it In order to classify  oriented rational knots, seen as numerator closures of
oriented rational tangles, we will always  compare compatible rational tangles.}  
\smallbreak

While the connectivity type of unoriented rational tangles may be $[0]$,  $[\infty]$
 or $[1],$ note that an oriented rational tangle of type I will have 
connectivity  type $[0]$ or
$[\infty]$ and an oriented rational tangle of type II will have  connectivity   type
$[0]$ or $[1].$ 

\bigbreak

$$ \picill4.6inby2.1in(R32) $$

\begin{center}
{ Figure 32 - Compatible and Incompatible Orientations } 
\end{center}
\vspace{3mm}


\noindent {\bf Bottom twist basics.} \ If two oriented tangles are incompatible,  adding
a  single half twist at the bottom of one of them yields  a new pair of compatible
tangles,  as Figure 32 illustrates. Note also that adding such a twist, although it
changes  the tangle, it does not change the isotopy type of the numerator closure. Thus,
up  to bottom twists, we are always able to compare oriented rational tangles of the same
orientation type. Further, note that if we add a positive bottom twist to an oriented 
rational tangle $T$ with fraction $F(T)= p/q$ we obtain the incompatible tangle
$T'=T*[+1]$ with fraction  $F(T') = 1/(1+ 1/F(T)) = p/(p+q).$
Similarly, if we add a negative twist we obtain the incompatible tangle $T''=T*[-1]$ with
fraction   $F(T'') = 1/(-1+ 1/F(T)) = p/(-p+q).$
It is worth noting here that the tangles $T'$ and $T''$ are compatible and $p+q \equiv
(-p+q) \, mod \, 2p,$ confirming the Oriented Schubert Theorem. 

\bigbreak

Schubert \cite{Sch2} proved his version of  Theorem \ref{Schubert2} 
 by using the $2$-bridge representation of rational knots and  links. We give a
tangle-theoretic proof of Schubert's Oriented Theorem, based upon the combinatorics of
the unoriented case and then analyzing how
 orientations and fractions are related. 

In our statement of Theorem \ref{Schubert2} in the introduction 
we restricted the denominators of the fractions to be odd.  This is a
restriction made for the 
purpose of comparison of tangles. There is no loss of generality, as will
be seen when we analyze the 
palindrome case in the proof at the end of this section. What happens is
this: In the case of $p$
odd and only one of $q$ and $q'$ even, one finds that the corresponding
tangles are incompatible.
We can then compare them by adding a bottom twist to one of the tangles.
 Adding this twist, the even denominator is replaced by an odd
denominator. In the case where $p$
is odd and both $q$ and $q'$ are even, one finds that the corresponding
tangles are compatible. 
In this case, we add a twist at the bottom of each tangle to preserve the
hypothesis that both 
denominators are odd. This extra twisting yields compatible tangles
and a successful comparison.

\smallbreak
The strategy of our proof is as follows. Consider an  
 oriented rational knot or link diagram $K$ given in standard form
 as $N(T),$ where $T$ is a rational tangle in continued fraction form.
  Our previous analysis tells us that, up to bottom twists, any other 
 rational tangle that closes to this knot is available as a cut on the
 given diagram. If two rational tangles close to give $K$ as an
 unoriented rational knot or link, then there are orientations on these tangles,
 induced from $K,$ so that the oriented tangles close
 to give $K$ as an oriented knot or link. Two tangles so produced may or may not be
compatible.  However, adding a bottom twist to one of two incompatible
 tangles results in two compatible tangles. {\it It is this possible twist
difference that gives rise to the change
 from modulus $p$ in the unoriented case to the modulus 
 $2p$ in the oriented case.}

We now analyze when, comparing with the
original standard cut, another cut
produces a compatible or incompatible tangle.  See Figure 34 for an example illustrating
the compatibility of orientations in the case of the palindrome cut. 
Note that reducing
all possible bottom twists implies $|p|>|q|$ for both tangles, if the two reduced 
tangles that we compare each time are compatible, or for only one, if they are
incompatible. Recall Figure 12 and the related analysis for the basic arithmetic of the
bottom twists.  
\bigbreak

\noindent {\bf Even bottom twists.} \   The simplest instance of the classification of
oriented rational knots is adding an
 {\it even number of twists} at the bottom of an oriented rational tangle $T.$ We then
 obtain a compatible tangle $T*1/[2n],$ and $N(T*1/[2n]) \sim N(T).$ If now $F(T) =
 p/q$, then  $F(T*1/[2n]) = F(1/([2n] + 1/T)) = 1/(2n+1/F(T))=  p/(2np + q).$
 Hence, if we set $2np + q = q'$ we have 

$$q\equiv q' \, mod \, 2p,$$ 
just as  Theorem \ref{Schubert2} predicts. 
 \bigbreak

  We then have to  compare the special cut and the palindrome cut with 
 the standard cut. Here also, the special cut is the easier to see whilst
 the palindrome cut requires a more sophisticated analysis. Figure 17 explained how to
 obtain the unoriented tangle of the special cut. Moreover, by Remark~2, adding a bottom
twist to the tangle of the special cut yields a tangle isotopic to the tangle of the
standard cut. 

  Figure 33 demonstrates that the special cut yields oriented
incompatible tangles.  More precisely, in the case of the special cut we are presented 
with the general fact that for any tangle $R$, $N([+1] + R)$ and $N([-1] -
1/R)$ are unoriented isotopic. With orientations coming from the 
cut we find that $S~=~[+1] + R$ and $S' =[-1] - 1/R$ are incompatible.
Adding a bottom twist yields oriented compatible tangles, which from the above are
isotopic. So, there is nothing to check and  the Oriented Schubert Theorem is
verified in the strongest possible way for the oriented special cut.

\bigbreak

$$ \picill5.5inby2.8in(R33) $$

\begin{center}
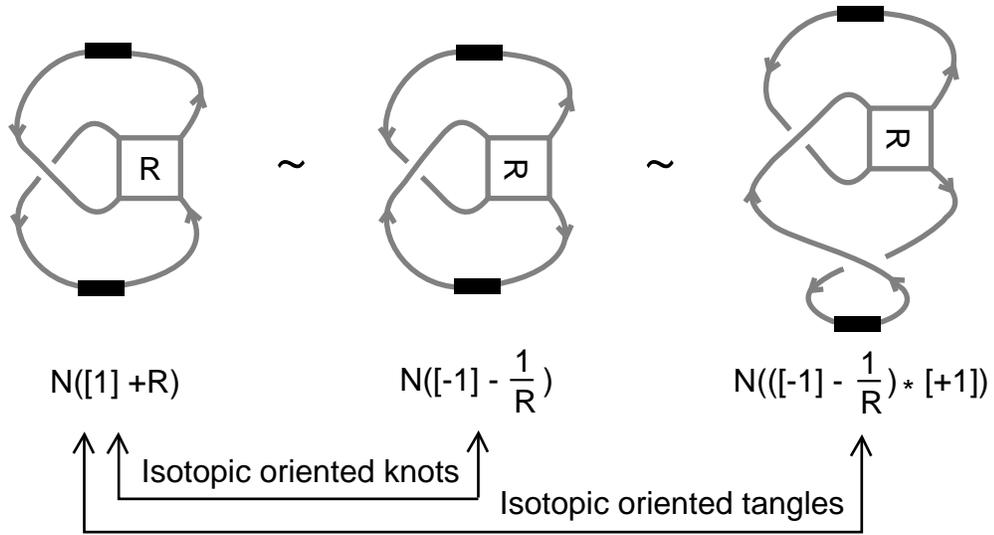
{ Figure 33 - The Oriented Special Cut yields Incompatible Tangles } 
\end{center}
\vspace{3mm}


We are left to examine the case of the palindrome cut. In order to
analyze this case, we must understand when the standard cut
and the palindrome cut are compatible or incompatible. Then we must
compare their respective fractions. Figure 34 illustrates how compatibility is
obtained by using a bottom twist, in the case of a palindrome cut.
 In this example we illustrate the standard and palindrome cuts on the
 oriented rational knot $K=N(T)=N(T')$ where
 $T=[[2],[1],[2]]$ and $T'$ its palindrome. As we can see, the two cuts place
 incompatible orientations on the tangles $T$ and
 $T'.$   Adding a twist at the bottom of $T'$ produces a tangle $T''=T'*[-1]$
 that is compatible with $T$. Now we compute
 $F(T) = F(T') = 8/3$ and $F(T'') = F(T'*[-1]) = 8/-5$ and we notice that  $3 \cdot (-5)
 \equiv 1 \, mod \, 16,$  as
 Theorem \ref{Schubert2} predicts. 

\bigbreak

$$ \picill5.4inby3in(R34) $$

\begin{center}
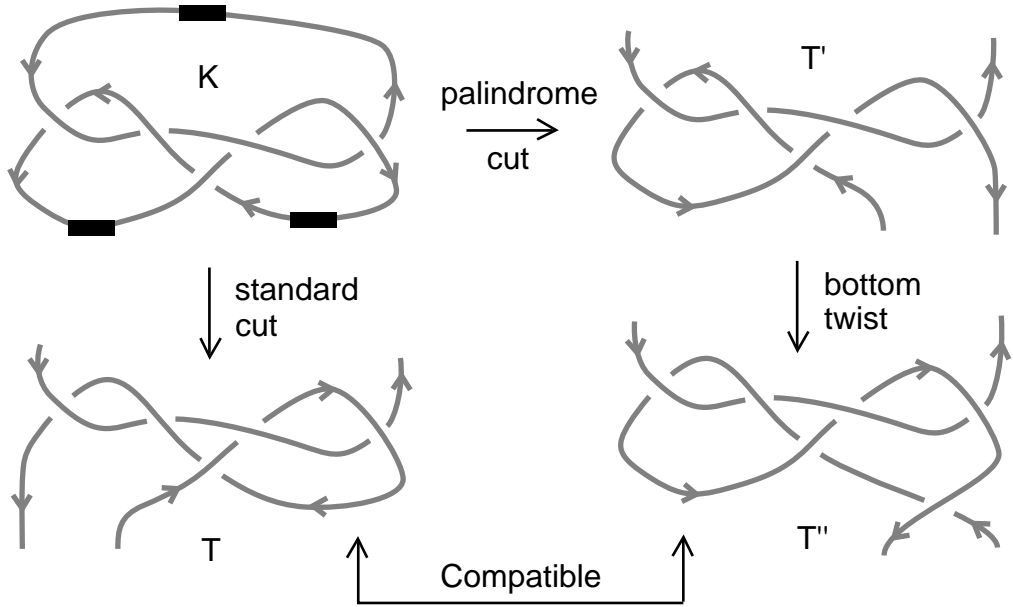
{ Figure 34 - Oriented Standard Cut and Palindrome Cut } 
\end{center}
\vspace{3mm}


 The study of the compatibility or not of the palindrome cut involves a deeper  
analysis  along the lines of Theorem 6. With the issues of connectivity in place we
can begin to analyze the different connectivities and parities in the 
standard and palindrome cuts on a rational knot or link in standard
$3$-strand-braid representation. See Figure 35.
In this figure we have enumerated the six connection structures for a
$3$-strand braid (corresponding to the six permutations 
of three points) with plat closures (of the braid augmented by an extra strand) corresponding to oriented
rational knots and links.  These closed connection patterns shall be called {\it connectivity charts.}  We
then show corresponding to each connectivity chart the  related standard and palindrome
cuts and the connectivity and parity of the corresponding tangles. Compatibility or
incompatibility of these tangles,  specified by an `i' or `c', can be read from the oriented 
diagrams in the figure. 

\bigbreak

$$ \picill6inby7.4in(R35) $$

\begin{center}
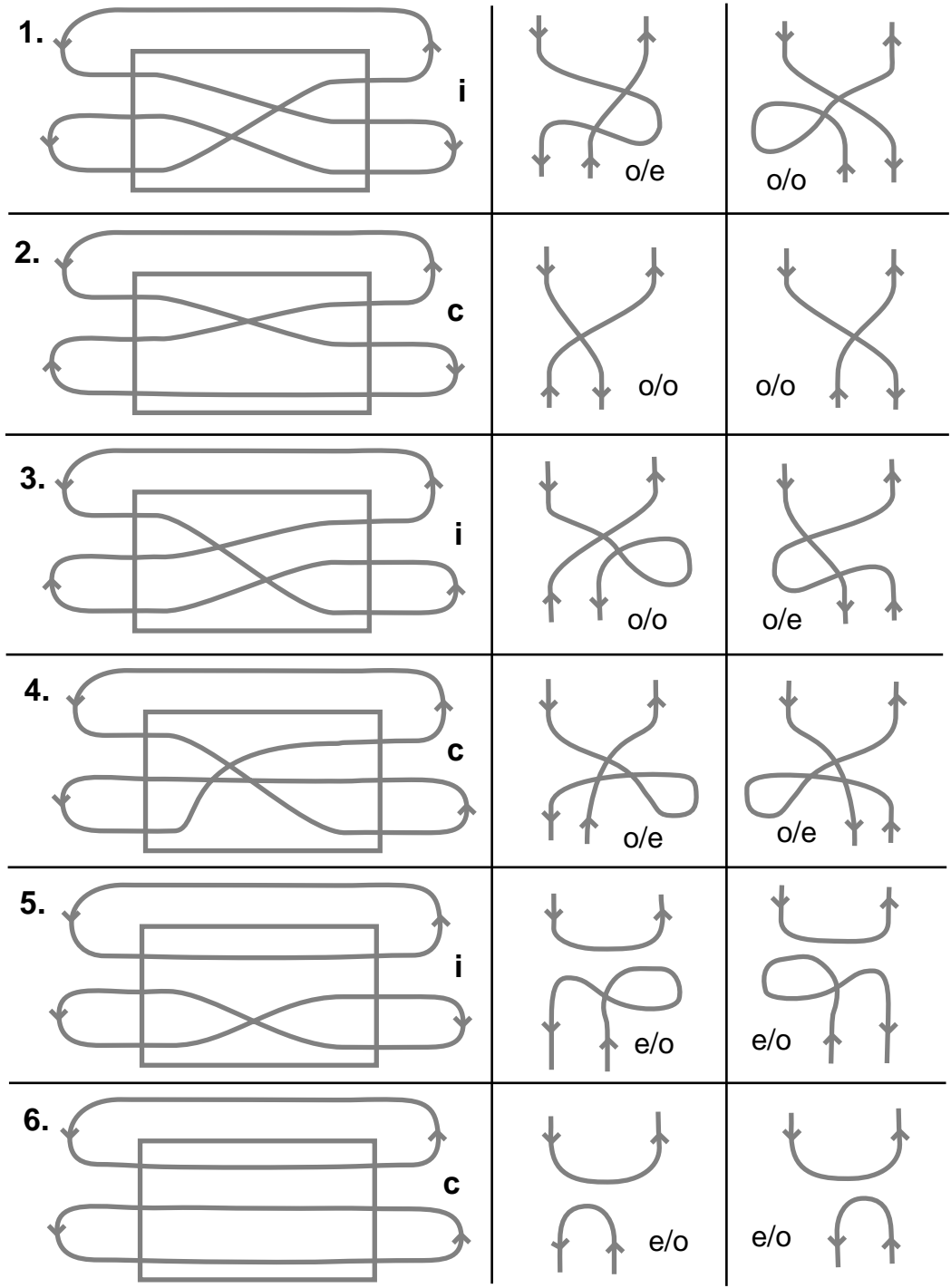
{ Figure 35 - The Six Connection Structures, Compatibility and Parity of the
Palindrome Cut } 
\end{center}
\vspace{3mm}


\noindent {\bf Proof of the Palindrome Cut.} It suffices to verify the Theorem in all
cases of the comparison of  standard and palindrome cuts on a rational knot $K$ in
continued fraction form. We can assume that $K=N([[a_{1}], \cdots ,[a_{n}]])$
with $n$ odd. Then  the tangle $T = [[a_{1}],\cdots,[a_{n}]]$ is,  by construction,  
the standard cut on $K.$ We know that the matrix product

$$M = M(a_{1})M(a_{2}) \cdots M(a_{n}) =
\left( \begin{array}{cc}
p & q' \\
q & u \\
\end{array} \right) $$

\noindent encodes the fractions of $T$ and its palindrome $T' = [[a_{n}], \cdots
,[a_{1}]],$ with $F(T) = p/q$ and $F(T') = p/q'.$  Note that, since
$Det(M) = -1,$ we have the formula  $$qq' = 1 + up$$ relating the
denominators of these fractions.
\bigbreak

\noindent {\bf Case 1. p odd, Part A:}
\smallbreak
\noindent {\it  If only one of $q$ or $q'$ is even} (parts 1 and 3 of
Figure 35), then the fact that 
$qq' = 1 + up$ implies the parity equation $e = 1 +uo$, hence {\em $u$ is
odd.}  Now refer to Figure 35 and note that
the standard and palindrome cuts are incompatible in both cases  1
and  3. (The cases are $\{
o/e,o/o \}$ and $\{ o/o, o/e \}.$) In order to obtain compatibility,
add a bottom twist to the cut with even denominator.  Without loss of 
generality, we may assume that $q'$ is even, so that we will compare $p/q$
and $p/(p + q').$ Note that 
$$q(p + q') = qp + qq' = qp + 1 + up = 1 + (u + q)p.$$
Since $u$ is odd and $q$ is odd, it follows that $(u + q)$ is even. Hence, 
$ q(p +q') \equiv 1 \, mod \, 2p,$ proving  Theorem \ref{Schubert2} in 
this case.
\bigbreak 

\noindent {\bf Case 1. p odd, Part B:}
\smallbreak
\noindent Now suppose that {\em both $q$ and $q'$ are even}. We are in
part  4 of Figure 35 and the two cuts are compatible.
Therefore we apply a bottom twist to each cut giving the fractions
$p/(p+q)$ and $p/(p+q')$ for comparison. Note that

$$(p+q)(p+q') = p^{2} + qp + q'p + qq' = 1 +(p + q + q' + u)p$$

\noindent and we have the parity equation

$$p + q + q' + u = o + e + e + o = e.$$

\noindent Hence $(p + q)(p +q') \equiv 1 \, mod \, 2p$ verifying the
Theorem in this case.
\bigbreak

\noindent {\bf Case 1. p odd, Part C:}
\smallbreak
\noindent Finally (for Case 1) suppose that {\em $q$ and $q'$ are both odd.}
Then the parity equation corresponding to \, 
$qq' = 1 + up$ \, is $$o = 1 + uo.$$ \noindent Hence $u$ is even so that $qq'
\equiv 1 \, mod \, 2p.$  We are in part  2 of Figure 35, 
and the standard and palindrome cuts are compatible. This is in accord with
the congruence above, hence the Theorem is verified in this 
case.
\bigbreak

\noindent {\bf Case 2. p even:}
\smallbreak
\noindent Now we assume that $p$ is even. This corresponds to parts 
5 and  6 in Figure 35 (two components). In part 5 the
cuts are compatible, while in part 6 the cuts are incompatible. In
either case, both $q$ and $q'$ are odd so that the fractions
$p/q$ and $p/q'$ both have the parity $e/o.$  The equation $qq' = 1 + up$
has corresponding parity equation $o = 1 + ue,$ and $u$ can be either
even or odd. In order to accomplish the proof of Case 2 we will show that 
\begin{enumerate}
\item $u$ is even if and only if the standard and palindrome cuts are
compatible.
\item $u$ is odd if and only if the standard and palindrome cuts are
incompatible.
\end{enumerate}
\smallbreak

\noindent We prove these statements by induction on the number of terms in
the continued fraction $[a_{1}, \cdots , a_{n}].$
The induction step consists in adding two more terms to the continued
fraction (thereby maintaining an odd number of terms).
That is, we shall examine a continued fraction in the form  $T_{n+2} =
[a_{1}, \cdots , a_{n+2}]$ that is {\em given} to be in cases  5 or 
 6 of Figure 35. See Figure 36. In Figure 36 the numbers that label the diagrams refer to
the cases in Figure 35.  We consider the structure of the ``predecessor"
$T_{n} = [a_{1}, \cdots , a_{n}]$  of $T_{n+2}$ which may be in the form 
 5  or  6, as shown in Figure 36 (in which case we can apply the induction
hypothesis) or it may be in one of the other four
cases  shown in Figure 36. 

\bigbreak

$$ \picill8inby5.4in(R36) $$

\begin{center}
{ Figure 36 - Inductive Connections } 
\end{center}
\vspace{3mm}


In Figure 36 we have shown the connectivity patterns
that result in $T_{n+2}$ landing in cases  5 or 6. In this 
figure the rectangular boxes indicate the internal connectivity of
$T_{n}$, and we have separated these specific cases into three types
labeled
$A$, $B$ and $C$ (not to be confused with subcases of this proof). 
In this Figure each case is labeled with the type of the predecessor. Thus in the 
$A$ row one sees the labels $3$ and $4$ because the boxed patterns are respectively of 
types $3$ and $4.$ In rows $A$ and $B$ the left hand entries are of type $6$ after the addition of
the two new terms, and the right hand entries  are of type $5.$ We then
check each of these cases to see that the induced value of $u$ in 
$T_{n+2}$ has the right parity. The calculations can be done by
multiplication of generating matrices for continued fractions just using the parity algebra. For
example, in Case A of Figure 36 we add two new odd parity terms to $T_{n}$ in order to obtain
$T_{n+2}.$ Thus we multiply the parity matrix for $T_{n}$ by the 
product 

$$\left( \begin{array}{cc}
o & 1\\
1 & 0\\
\end{array} \right)
\left( \begin{array}{cc}
o & 1\\
1& 0\\
\end{array} \right)
=
\left( \begin{array}{cc}
e & o\\
o& 1\\
\end{array} \right)$$

\noindent in order to obtain the parity matrix for $T_{n+2}.$ 
\smallbreak

\noindent In particular, if $T_{n}$ is in case  3 of Figure 35,
then it has fraction parities $o/o$ and $o/e$ and hence parity
matrix
$\left( \begin{array}{cc}
o & e\\
o & o\\
\end{array} \right).$  Multiplying this by 
$\left( \begin{array}{cc}
e & o\\
o & 1\\
\end{array} \right),$ we obtain 
$$\left( \begin{array}{cc}
o & e\\
o & o\\
\end{array} \right)
\left( \begin{array}{cc}
e & o\\
o & 1\\
\end{array} \right)
=
\left( \begin{array}{cc}
e & o\\
o & e\\
\end{array} \right).$$
\noindent Thus the new $u$ for $T_{n+2}$ is even. Since the connectivity
diagram for $T_{n+2}$ in this case, as shown in Figure 36,
has compatible standard and palindrome cuts, this result for the parity of
$u$ is one step in the verification of the induction 
hypothesis. Each of the six cases is handled in this same way. We omit the
remaining details and assert that the values of $u$ obtained 
in each case are correct with respect to the connection structure. This
completes the proof of Case 2.
\smallbreak

Since Cases 1 and 2 encompass all the different possibilities for the
standard and palindrome cuts, this completes the proof of 
the Oriented Schubert Theorem. $\hfill \Box $

 \section{Strongly Invertible Links}

 An oriented knot or link is invertible if it is oriented
 isotopic to the link obtained from it by reversing the orientation of each component.
We have seen (Lemma \ref{invertible}) that  rational knots and
links are invertible.  A link $L$ of two components is said to
be {\it strongly invertible} if $L$ is ambient isotopic to itself with the orientation
of only one component reversed. In Figure 37 we illustrate the link
$L=N([[2],[1],[2]]).$ This is a strongly invertible link as is apparent by a $180^0$
vertical rotation.  This link is well-known as the Whitehead link, a link with linking
number zero. Note that since
$[[2], [1], [2]]$ has fraction equal to $1 + 1/(1 + 1/2) = 8/3$ this link is non-trivial
via the classification of rational knots and links. Note also that $3 \cdot 3 = 1 + 1
\cdot 8.$ 

\bigbreak

$$ \picill4.3inby2in(R37) $$

\begin{center}
{ Figure 37 - The Whitehead Link is Strongly Invertible } 
\end{center}
\vspace{3mm}

 
 In general we have the following.

 \begin{th}\label{strongly} \ Let $L=N(T)$ be an oriented rational link with associated
 tangle fraction $F(T) = p/q$  of parity $e/o,$ with $p$ and $q$ relatively prime and
$|p|>|q|.$ Then $L$ is strongly invertible if and only if $q^{2} = 1 + up$ with $u$  an
odd integer. It follows that strongly invertible links are all numerators of rational
tangles  of the form $[[a_1], [a_2], \ldots, [a_k], [\alpha], [a_k], \ldots, [a_2],
[a_1]]$ for any integers $a_1, \ldots, a_k, \alpha.$
  \end{th}

 \noindent {\em Proof.}  In $T$ the upper two end arcs close to form one
 component of $L$ and the lower two end arcs close to form the other  component of $L.$
Let $T'$ denote the tangle obtained from the oriented tangle $T$ by reversing the
orientation of the component containing the lower two
 arcs and let $N(T') = L'.$ (If $T''$ denotes the tangle obtained from the oriented
tangle $T$ by reversing the orientation of the component containing the upper two
 arcs we have seen that by a vertical $180^0$ rotation the link $N(T')$ is isotopic to 
the link $N(T'').$ So, for proving Theorem \ref{strongly} it suffices to consider only the
case above.) 

 Note that $T$ and $T'$ are incompatible.
  Thus to apply Theorem \ref{Schubert2} we need to perform a bottom twist on
$T'.$ Since $T$ and $T'$ have the same fraction $p/q,$ after applying the twist we need
to compare the  fractions $p/q$ and $p/(p+q).$ Since $q$ is
 not congruent to $(p+q)$ modulo $2p$, we need to determine when $q(p+q)$ is
 congruent to $1$ modulo $2p.$ This will happen exactly  when $qp + q^{2} = 1 + 2Kp$ for
some integer $K.$ The last equation is the same as saying that $q^{2} = 1 + up$ with $u
= 2K-q$ odd, since $q$ is odd. Now it follows from the Palindrome Theorem for continued fractions
that 
$q^{2} = 1 + up$ with $u$ odd and $p$ even if and only if the fraction $p/q$ with
$|p|>|q|$ has a palindromic  continued  fraction expansion with an odd  number of terms
(the proof is the same in form as the corresponding argument given in the proof of
Theorem $5$). 
That is, it has a continued fraction
in the form
$$[a_{1}, a_{2}, \cdots, a_{n}, \alpha, a_{n}, a_{n-1}, \cdots, a_{2}, a_{1}].$$
It is then easy to see that the corresponding rational link  is ambient isotopic to
itself through a vertical $180^0$ rotation. Hence it is strongly invertible. 
It follows from this that all strongly invertible rational links are
ambient isotopic to themselves through a $180^0$ rotation just as
in the example of the Whitehead link given above. 
 This  completes the proof  of the Theorem. $\hfill \Box$ 

\begin{rem}{\rm \  Excluding the possibility $T = [\infty],$ as  $F(T) =1/0$ does not
have the parity $e/o,$ we may assume $q \neq 0.$ And since $q$ is odd (in order that the rational
tangle has two components), the integer
$u = 2K-q$ in the equation  $q^{2} = 1 + up$ cannot be zero. It follows then that the links of
the type $N([2n]),$ for $n\in \ZZ, \ n\neq 0$ with tangle fraction $2n/1$ are not
invertible (recall the example in  Figure 30). Note that, for $n=0$ we have
$T = [0]$ and  $F(T) =0/1,$ and in this case Theorem \ref{strongly} is confirmed, since 
 $1^{2} = 1 + u\,0,$ for any $u$ odd.  See Figure 38 for another example of a strongly
invertible link.  In this case the link is $L = N([[3], [1], [1], [1], [3]])$ with $F(L)
= 40/11.$ Note that
$11^{2} =1 + 3\cdot 40,$ fitting the conclusion of Theorem \ref{strongly}.
} \end{rem}

\bigbreak

$$ \picill2.5inby1.6in(R38) $$

\begin{center}
{ Figure 38 - An Example of a Strongly Invertible Link } 
\end{center}
\vspace{3mm}


 \noindent {\bf Acknowledgments.} \ A part of this work and the work in \cite{KL1} was
done at  G\"ottingen Universit\"at. Both authors
acknowledge with gratitude the research facilities offered there. 
It also gives us  pleasure to acknowledge a list of
places and meetings where we worked on these matters. These are: G\"{o}ttingen, Santa Fe, Chicago,
Austin, Toulouse, Korea, Athens,  Marne-la-Vall\'{e}e, Siegen, Las Vegas, Marseille,
Berkeley, Dresden, Potsdam NY.


 \bigbreak

 \noindent {\sc L.H.Kauffman: Department of Mathematics, Statistics and
 Computer Science, University
 of Illinois at Chicago, 851 South Morgan St., Chicago IL 60607-7045,
 U.S.A.}

 \vspace{.1in}
 \noindent {\sc S.Lambropoulou: Department of Mathematics, 
 National Technical University of Athens,
 Zografou campus, GR-157 80 Athens, Greece. }

\vspace{.1in}
\noindent {\sc E-mails:} \ {\tt kauffman@math.uic.edu  \ \ \ \ \ \ \ \ sofia@math.ntua.gr 

\noindent http://www.math.uic.edu/$\tilde{~}$kauffman/ \ \ \ \ http://users.ntua.gr/sofial}



 \end{document}